\documentclass[10pt,a4paper]{article}
\usepackage{amsfonts,amssymb,amsthm,amsmath,amsopn}
\usepackage{geometry}
\usepackage{hyperref}
\usepackage[monochrome]{xcolor}
\usepackage{graphicx}
\usepackage{subfigure}

\theoremstyle{remark}
	\newtheorem*{remark}{Remark}
\theoremstyle{definition}
	\newtheorem{definition}{Definition}[section]
\theoremstyle{definition}

\title{
	Geometric Optimal Trajectory Tracking of Nonholonomic Mechanical Systems
	}

\author{
	Leonardo J. Colombo
	\thanks{ L. Colombo and D. Mart\'in de Diego are members of the Instituto de Ciencias Matem\'aticas (CSIC-UAM-UC3M-UCM), Calle Nicol\'as Cabrera 15, Campus UAM, Cantoblanco, 28049, Madrid, Spain. (\email{leo.colombo@icmat.es, david.martin@icmat.es}).}
	\and
	David Mart\'in de Diego\footnotemark[1]
	\and
	Aradhana Nayak
	\thanks{A. Nayak is a member of the Centre Automatique et Syst\'em\'es (MINES ParisTech). 60 Boulevard Saint-Michel, 75006 Paris, France. (\email{aradhana.nayak@mines-paristech.fr}).}
	\and
	Rodrigo T. Sato Mart\'in de Almagro \thanks{ R. T. Sato Mart\'in de Almagro is a member of the Lehrstuhl f{\"u}r Technische Dynamik (FAU), Immerwahrstr. 1, 91058 Erlangen, Germany (\email{rodrigo.t.sato@fau.de})}
	}

\newcommand{\R}{\mathbb{R}}
\newcommand{\Flder}{\rightarrow}

\newcommand{\ra}{\rightarrow}

\newcommand{\lcf}{\lbrack\! \lbrack}
\newcommand{\rcf}{\rbrack\! \rbrack}

\newcommand{\proa}{A^*G \mbox{$\;$}_{\tau^*} \kern-3pt\times_\alpha
G \mbox{$\;$}_\beta \kern-3pt\times_{\tau^*} A^*G}

\newcommand{\email}[1]{\protect\href{mailto:#1}{#1}}

\newcommand\keywordsname{Key words}
\newcommand\keywordname{Key word}
\newcommand\AMSname{AMS subject classifications}
\newcommand\AMname{AMS subject classification}

\newenvironment{@abssec}[1]{%
     \if@twocolumn
       \section*{#1}%
     \else
       \vspace{.05in}\footnotesize
       \parindent .2in
         {\upshape\bfseries #1. }\ignorespaces 
     \fi}
     {\if@twocolumn\else\par\vspace{.1in}\fi}

\newenvironment{keywords}{\begin{@abssec}{\keywordsname}}{\end{@abssec}}

\newenvironment{AMS}{\begin{@abssec}{\AMSname}}{\end{@abssec}}

\begin{document}

\maketitle

\begin{abstract}
We study the tracking of a trajectory for a nonholonomic system by recasting the problem as a constrained optimal control problem. The cost function is chosen to minimize the error in positions and velocities between the trajectory of a nonholonomic system and the desired reference trajectory, both evolving on the distribution which defines the nonholonomic constraints.  The problem is studied from a geometric framework. Optimality conditions are determined by the Pontryagin Maximum Principle and also from a variational point of view, which allows the construction of geometric integrators. Examples and numerical simulations are shown to validate the results.\end{abstract}

\begin{keywords}
  Optimal control, Trajectory planning, Nonholonomic systems, Variational integrators.
\end{keywords}

\begin{AMS}
  22E70, 37K05, 37J15, 37M15, 37N35, 49J15, 91B69, 93C10.
\end{AMS}

\section{Introduction}

Nonholonomic optimal control problems arise in many engineering applications,
for instance systems with wheels, such as cars and bicycles, and systems with
blades or skates. There are thus multiple applications in the context of wheeled
motion, space or mobile robotics and robotic manipulation. The earliest work on control of nonholonomic systems is by R. W. Brockett in \cite{brocket}. A. M. Bloch \cite{Bl}, \cite{bl1} has examined several control theoretic issues which pertain to both holonomic and nonholonomic systems in a very general form. The seminal works about stabilization in nonholonomic control systems were done by A. M. Bloch, N. H. McClamroch, and M. Reyhanoglu in \cite{bl1}, \cite{bl2}, \cite{bl3}, \cite{blochnonh}, and more recent results on the topic has been developed by A. Zuyev \cite{zh}.

    Geometrically, a conservative dynamical system of mechanical type is completely determined by a Riemannian manifold $Q$, the kinetic energy of the mechanical system, which is defined through the Riemannian metric $\mathcal{G}$ on $Q$ and the potential forces encoded into a potential (conservative) function $V:Q\to\mathbb{R}$. These objects, together with a non-integrable distribution $\mathcal{D}\subset TQ$ on the tangent bundle of the configuration space determines a nonholonomic mechanical system (see \cite{Bl} and references therein). Note that the description that we propose for dynamical systems of mechanical type only  apply for conservative systems, as there might be also non-conservative (dissipative and gyroscopic) forces in general mechanical systems.

Stabilization of an equilibrium point of a mechanical system on a Riemannian manifold has been a problem well studied in the literature from a geometric framework along the last decades (see \cite{Bl} and \cite{bullolewis} for a review on the topic). Further extensions of these results to the problem of tracking a smooth and bounded trajectory can be
found in \cite{bullolewis} where a proportional and derivative plus feed forward
(PD+FF) feedback control law is proposed for tracking a trajectory on a Riemannian manifold using error functions. 

For trajectory tracking, the usual approach of stabilization of error dynamics \cite{kodi}, \cite{aradhana1}, \cite{aradhana2}, \cite{sanyal}  cannot be utilized for nonholonomic systems. This is because there does not exist a $\mathcal{C}^1$ (even continuous) state feedback which can stabilize the trajectory of a nonholonomic system about a desired equilibrium point. The closed loop trajectory violates Brockett's condition \cite{brocket2}, \cite{blochnonh} which states that any system of the form $\dot{x} = f(x,u)$ must have a neighborhood of zero in the image of the map $x \to f(x,u)$ for some $u$ in the control set.  This result appears in Theorem 4 in \cite{blochnonh}.

In this paper, we introduce a geometrical framework in nonholonomic mechanics to study tracking of trajectories for nonholonomic systems based on \cite{leothesis}, \cite{Cobook}, \cite{CoMa}. The application of modern tools from differential geometry in the fields of mechanics, control theory and numerical integration has led to significant progress in these research areas. For instance, the study on the geometrical formulation of the nonholonomic equations of motion has led to better understanding of different engineering problems such locomotion generation, controllability, motion planning, and trajectory tracking \cite{Bl},  \cite{bullolewis}, \cite{Cobook}. 

Combining the ideas of geometric methods in control theory, nonholonomic systems and optimization techniques, in this paper, we study the underlying geometry of a tracking problem for nonholonomic systems by understanding it as a constrained optimal control problem for mechanical systems subject to nonholonomic constraints.

Given a reference trajectory $\gamma_{r}(t)=(q_r(t), v_r(t))$ on $\mathcal{D}$ the problem studied in this work consists on finding an admissible curve $\gamma(t)\in\mathcal{D}$, solving a dynamical control system, with prescribed boundary conditions on $\mathcal{D}$ and minimizing a cost functional which involves the error  between the reference trajectory and the trajectory one wants to find (in terms of both, positions and velocities), and the effort of the control inputs. This cost functional is accomplished with a weighted  terminal cost  (also known as Mayer term) which induces a constraint into the dynamics on $\mathcal{D}$.  %

We propose a geometric derivation of the equations of motion for tracking a trajectory of a nonholonomic system as an optimal
control problem from two different points of view: as a constrained optimal control problem on the tangent space to the distribution $\mathcal{D}$ and from the Pontryagin Maximum Principle (PMP), where the optimal Hamiltonian is defined on the cotangent bundle of the constraint distribution. Both approaches allow the reduction in the degrees of freedom of the equations for the optimal control problem, compared with typical methods describing the dynamics of a nonholonomic system, as the ones arising from  the application of the classical Lagrange-d'Alembert principle. The main advantages in this geometric framework consist in the use of a basis of vector fields adapted to $\mathcal{D}$ allowing such a reduction of some degrees of freedom in the dynamics for a nonholonomic mechanical system.

 It is well known that (see \cite{Bl} for instance) Hamilton equations (in the cotangent bundle), are the dual representation of Euler-Lagrange equations (in the tangent bundle). By employing an arbitrary discretization of the necessary conditions for optimality arising from the PMP together with a shooting method for the boundary value problem, one can observe that for mechanical systems, the physical behavior of the system is not respected. Therefore it is needed to develop numerical algorithms showing a good qualitative behavior of solutions in simulations.  Our motivation to develop a Lagrangian formalism for the optimal trajectory tracking problems is mainly based on the fact that by considering a Lagrangian formalism it is possible to construct variational integrators. That is, a class of geometric numerical schemes that preserves the qualitative features of the system such as momentum preservation and symplecticity, and have remarkably good long-time energy behavior. This can be achieved by discretizing the variational principle, instead of discretizing the equations of motion as is usual in the literature to construct numerical methods for this class of problems. Moreover, it is also well known that Noether's theorem (given in the Lagrangian framework) provides a direct link between symmetries and conserved quantities which is preserved by the discretization of variational principles in the Lagrangian framework.

To test the efficiency of the proposed approach with the PMP, we use a Runge Kutta integrator together with a shooting method in the solution of a trajectory optimization for a simple but challenging benchmark mechanical system: a fully actuated particle subject to a nonholonomic constraint into the dynamics. We observed in the simulations how difficult is to achieve the reference trajectory in the constraint submanifold under the boundary conditions in the problem set-up. This motivate to us to propose a new numerical scheme to achieve the reference trajectory. This new scheme is based on a variational integrator. Such an integrator is tested in a classical nonholonomic system of mechanical type: the Chaplyigin sleigh. Numerical simulations exhibit an accurate convergence to the reference trajectory and a good behavior of the energy associated with the optimal control problem. Preliminaries results of this work by employing the PMP can be found our conference paper \cite{ACC}.

The paper is structured as follows: we introduce mechanical systems on a manifold, connections on a Riemannian manifold and the geometry of nonholonomic dynamical systems on Section \ref{secprel}, together with the examples we used as benchmarks: the nonholonomic particle and the Chaplygin sleigh. Section \ref{section3} introduces the details of the problem under study motivated by the non-existence of a $\mathcal{C}^1$ feedback control to stabilize the error dynamics in nonholonomic systems. Necessary conditions for extrema in the proposed optimal control problem are studied from the PMP and from a variational formalism in Section \ref{sec4}. The last motivate the construction of variational integrators in Section \ref{sec5}. We also show numerical results and analyze the results we obtain.  A final discussion and further applications and extensions of this work are presented in Section \ref{sec6} 

\section{Nonholonomic mechanical systems}\label{secprel}
Let $Q$ be the configuration space of a mechanical system, a differentiable manifold with $\dim(Q)=n$, and local coordinates denoted by $(q^i)$ for $i=1,\ldots,n$. Most nonholonomic systems have linear constraints on velocities, and these are the ones we will consider. Linear constraints on the
velocities (or Pfaffian constraints) are locally given by equations of the form $\phi^{a}(q^i, \dot{q}^i)=\mu^a_i(q)\dot{q}^i=0, \, 1\leq a\leq
m$, depending, in general, on their configurations and  their
velocities. 

From an intrinsic point of view, the linear
constraints are defined by a regular distribution ${\mathcal D}$ on
$Q$ of constant rank $(n-m)$ such that the annihilator of ${\mathcal
D}$ is locally given at each point of $Q$ by
${\mathcal D}^o_{q} = \operatorname{span}\left\{ \mu^{a}(q)=\mu_i^{a}dq^i \; ; 1 \leq a
\leq m \right\}$, where $\mu^{a}$ are independent one-forms at each point of $Q$.

We restrict ourselves to the case of nonholonomic mechanical systems where the Lagrangian is of mechanical type, that is, mechanical systems with a dynamics described by a Lagrangian function $L:TQ\to\R$ which is defined by \[
L(v_q)=\frac{1}{2}\mathcal{G}(v_q, v_q) - V(q),
\]
with $v_q\in T_qQ$, where $\mathcal{G}$ denotes a Riemannian metric on $Q$ representing the kinetic energy of the systems, and
$V:Q\ra\R$ is a potential function. 

Assume that the Lagrangian system is subject to nonholonomic constraints, defined by a regular distribution $\mathcal{D}$ on $Q$ with corank$(\mathcal{D})=m$.
Denote by $\tau_{\mathcal{D}}:\mathcal{D}\ra Q$ the canonical
projection from $\mathcal{D}$ to $Q$, denote by 
$\Gamma(\tau_{\mathcal{D}})$ the set of sections of $\tau_{D}$ and also denote by $\mathfrak{X}(Q)$ the set of vector fields taking values on $\mathcal{D}.$ If $X, Y\in\mathfrak{X}(Q),$ then
$[X,Y]$ denotes the standard Lie bracket of vector fields.

\begin{definition}\label{nonholonomicsystem}
A \textit{nonholonomic mechanical system} on a smooth manifold $Q$ is given
by the triple $(\mathcal{G}, V, \mathcal{D})$, where $\mathcal{G}$ is
a Riemannian metric on $Q,$ representing the kinetic energy of the
system, $V:Q\ra\R$ is a smooth function representing the potential
energy and $\mathcal{D}$ a non-integrable smooth distribution on $Q$
representing the nonholonomic constraints.
\end{definition}

Given $X,Y\in\Gamma(\tau_{\mathcal{D}})$ that is,
$X(x)\in\mathcal{D}_{x}$ and $Y(x)\in\mathcal{D}_{x}$ for all $x\in
Q,$ then it could happen that $[X,Y]\notin\Gamma(\tau_{\mathcal{D}})$
since $\mathcal{D}$ is nonintegrable. We want to obtain a bracket definition for sections on $\mathcal{D}.$
Using the Riemannian metric $\mathcal{G}$ we can define two
complementary orthogonal projectors ${\mathcal P}\colon TQ\to {\mathcal D}$ and ${\mathcal Q}\colon TQ\to {\mathcal
D}^{\perp},$ with respect to the tangent bundle orthogonal decomposition $\mathcal{D}\oplus\mathcal{D}^{\perp}=TQ$. Therefore, given $X,Y\in\Gamma(\tau_{\mathcal{D}})$ we define the
\textit{nonholonomic bracket}
$\lcf\cdot,\cdot\rcf:\Gamma(\tau_{\mathcal{D}})\times\Gamma(\tau_{\mathcal{D}})\rightarrow\Gamma(\tau_{\mathcal{D}})$
as $\lcf X_A,X_B\rcf:=\mathcal{P}[X_A,X_B]$. This Lie bracket verifies the usual properties of a Lie bracket except the Jacobi identity (see \cite{BlCoGuMdD}, \cite{leo2} for example). 

\begin{definition}
Consider the restriction of the Riemannian metric $\mathcal{G}$ to
the distribution $\mathcal{D}$, $\mathcal{G}^{\mathcal{D}}:\mathcal{D}\times_{Q}\mathcal{D}\ra\R$
and define 
$\displaystyle{\nabla^{\mathcal{G}^{\mathcal{D}}}:\Gamma(\tau_{\mathcal{D}})\times\Gamma(\tau_{\mathcal{D}})\ra\Gamma(\tau_{\mathcal{D}})}$, the \textit{Levi-Civita connection}
determined by the following two properties:

\begin{enumerate}
\item $\lcf X,Y\rcf=\nabla_{X}^{\mathcal{G}^{\mathcal{D}}}Y-\nabla_{Y}^{\mathcal{G}^{\mathcal{D}}}X,$
\item $X(\mathcal{G}^{\mathcal{D}}(Y,Z))=\mathcal{G}^{\mathcal{D}}(\nabla_{X}^{\mathcal{G}^{\mathcal{D}}}Y,Z)+\mathcal{G}^{\mathcal{D}}(Y,\nabla_{X}^{\mathcal{G}^{\mathcal{D}}}Z).$
\end{enumerate}

\end{definition}

Let $(q^{i})$ be local coordinates on $Q$ and $\{e_{A}\}$ be independent vector fields on
$\Gamma(\tau_{D})$ (that is, $e_A(x)\in {\mathcal D}_x$) such that $\mathcal{D}_{x}=\hbox{span }\{e_{A}(x)\},
\, x\in U\subset Q$. Then, we can determine the \textit{Christoffel
symbols} $\Gamma_{BC}^{A}$ associated with the connection
$\nabla^{\mathcal{G}^{\mathcal{D}}}$ by $\displaystyle{\nabla_{e_{B}}^{\mathcal{G}^{\mathcal{D}}}e_{C}=\Gamma_{BC}^{A}(q)e_{A}.}$ Note that the coefficients $\Gamma_{AB}^{C}$ of the connection  $\nabla^{{\mathcal
G}^{\mathcal{D}}}$ can be also computed by (see \cite{CoMa} for details) \begin{equation}\label{relation}\Gamma_{AB}^{C}=\frac{1}{2}(\mathcal{C}_{CA}^{B}+\mathcal{C}_{CB}^{A}+\mathcal{C}_{AB}^{C})\end{equation}where the constant structures $\mathcal{C}_{AB}^{C}$ are defined by $\lcf X_A,X_B\rcf=\mathcal{C}_{AB}^{C}X_C$.

\begin{definition}
A curve $\gamma:I\subset\R\ra\mathcal{D}$ is \textit{admissible} if
$\displaystyle{\gamma(t)=\frac{d\sigma}{dt}(t)}$, where $\tau_{\mathcal{D}}\circ\gamma=\sigma$.
\end{definition}

Given local coordinates on $Q,$ $(q^{i})$ with $i=1,\ldots,n;$ and
$\{e_{A}\}$ sections on $\Gamma(\tau_{\mathcal{D}})$, with $A=1,\ldots,n-m$, such that
$\displaystyle{e_{A}=\rho_{A}^{i}(q)\frac{\partial}{\partial
q^{i}}}$ we introduce induced coordinates $(q^{i},v^{A})$ on
$\mathcal{D}$, where, if $e\in\mathcal{D}_{x}$ then
$e=v^{A}e_{A}(x).$ Therefore, the curve $\gamma(t)=(q^{i}(t),v^{A}(t))$ is
admissible if $\dot{q}^{i}(t)=\rho_{A}^{i}(q(t))v^{A}(t)$.

Consider the restricted Lagrangian function
$\ell:\mathcal{D}\rightarrow\mathbb{R},$
$$\ell(v)=\frac{1}{2}\mathcal{G}^{\mathcal{D}}(v,v)-V(\tau_{D}(v)),\hbox{ with }  v\in\mathcal{D}.$$

\begin{definition}
A \textit{solution of the nonholonomic problem} is an admissible
curve $\gamma:I\rightarrow\mathcal{D}$ such that
$$\nabla_{\gamma(t)}^{\mathcal{G}^{\mathcal{D}}}\gamma(t)+grad_{\mathcal{G}^{\mathcal{D}}}V(\tau_{\mathcal{D}}(\gamma(t)))=0.$$
\end{definition} Here the section $grad_{{\mathcal G}^{\mathcal{D}}}V\in\Gamma(\tau_{\mathcal{D}})$ is characterized by \[ {{\mathcal G}^{\mathcal{D}}}(grad_{{\mathcal
G}^{\mathcal{D}}}V, X) = X(V), \; \mbox{ for  every } X \in
\Gamma(\tau_{\mathcal{D}}).
\]

These equations are equivalent to the \textit{nonholonomic
equations}. Locally, these equations are given by
\begin{align}
\dot{q}^{i}&=\rho_{A}^{i}(q)v^{A}\label{eqq1}\\
\dot{v}^{C}&=-\Gamma_{AB}^{C}v^{A}v^{B}-(\mathcal{G}^{\mathcal{D}})^{CB}\rho_{B}^{i}(q)\frac{\partial V}{\partial q^{i}},\label{eqq2}
\end{align} where $(\mathcal{G}^{\mathcal{D}})^{AB}$ denotes the coefficients of the inverse matrix of $(\mathcal{G}^{\mathcal{D}})_{AB}$ determined by $\mathcal{G}^{\mathcal{D}}(e_{A},e_{B})=(\mathcal{G}^{\mathcal{D}})_{AB}.$

\begin{remark}
The nonholonomic  equations \eqref{eqq1}-\eqref{eqq2} only depend on the coordinates
$(q^{i},v^{A})$ on $\mathcal{D}$.  Therefore the nonholonomic
equations are free of Lagrange multipliers. These equations are
equivalent to the  \textit{nonholonomic Hamel equations} (see
\cite{BloZen}, for example).

\end{remark}

\subsection{Example: the Chaplygin sleigh }\label{chaplygin}

The Chaplygin sleigh (see \cite{Bl})  is a rigid body moving on a horizontal plane with three contact points, two of which slide freely without friction. The third one is a knife edge, which imposes the nonholonomic constraint of no motion perpendicular to the direction of the blade. The configuration space is $Q=SE(2)$, with local coordinates $(x_1,x_2,\theta)$. The coordinates $(x_1,x_2)$ denote the contact point of the blade with the plane and $\theta$ the orientation of the blade.

\begin{figure}[h!]
  \begin{center}\includegraphics[width=6.7cm]{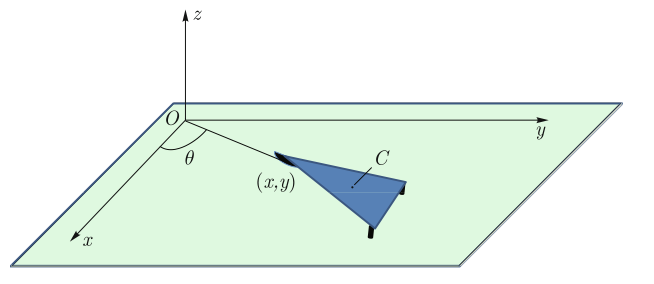}
  \caption{The Chaplygin sleigh}\end{center}\end{figure}

The Lagrangian is of kinetic type and  if we assume that the center of mass lies in the line through the blade then it is given by
$$
L=\frac{1}{2}\left( (J+ma^2)\dot{\theta}^2+m\left(\dot{x}_1^2+\dot{x}_2^2+2a\dot{\theta}(-\dot{x}_1\sin\theta+\dot{x}_2\cos\theta)\right) \right),
$$
where $m$ denotes the mass of the body, $J$ the moment of inertia relative to the center of mass and $a$ the distance between the center of mass and the contact point of the blade. The matrix of the metric defining the kinetic Lagrangian is given by
$$
\left(
\begin{array}{ccc}
m & 0 & -ma\sin\theta \\
0 & m & ma\cos\theta \\
-ma\sin\theta & ma\cos\theta & J+ma^2
\end{array}
\right) \, .
$$
The nonholonomic constraint is $\dot{x}_2\cos(\theta)=-\dot{x}_1\sin(\theta)$, which defines a non-integrable distribution 
$$
{\mathcal D}=\mbox{span}\left\{\frac{\partial}{\partial \theta},  \cos\theta\frac{\partial}{\partial x_1}+\sin\theta\frac{\partial}{\partial x_2}  \right\} .
$$

To derive the nonholonomic equations in adapted coordinates, we choose the following orthonormal basis adapted to $\mathcal{D}$ and $\mathcal{D}^{\bot}$ :
\begin{align*}
{\mathcal D}&=\mbox{span}\left\{X_1=\frac{1}{\sqrt{J+ma^2}}\frac{\partial}{\partial \theta}\, , \, X_2=\frac{1}{\sqrt{m}}\left( \cos\theta\frac{\partial}{\partial x_1}+\sin\theta\frac{\partial}{\partial x_2} \right) \right\} \, ,\\
{\mathcal D}^{\bot}&=\mbox{span}\left\{ X_3=\Gamma\left(\frac{(J+ma^2)}{ma}\sin\theta\frac{\partial}{\partial x_1}-\frac{(J+ma^2)}{ma}\cos\theta\frac{\partial}{\partial x_2}+\frac{\partial}{\partial \theta}\right) \right\} \, ,
\end{align*}with $\Gamma=\frac{1}{\sqrt{\frac{(J+ma^2)^2}{ma^2}-(J+ma^2)}}$. Denote by $(q^i,v^A)=(x_1,x_2,\theta,v^1,v^2,v^3)$ the induced coordinates. A straightforward computation shows that the functions $\rho_A^i$ are given by $$\rho_1^1=\rho_1^2=\rho_2^3=0,\quad \rho_1^3=\frac{1}{\sqrt{J+ma^2}},\quad\rho_2^1=\frac{\cos\theta}{\sqrt{m}},\quad\rho_2^2=\frac{\sin\theta}{\sqrt{m}}.$$

In the induced coordinates, the restricted Lagrangian $\ell:{\mathcal D}\longrightarrow \mathbb{R}$ is given by $\ell(q^i,v^A)=\frac{1}{2}((v^1)^2+(v^2)^2)$, and the nonholonomic constraint by $v^3=0$.

The the nonholonomic equations for the Chaplygin sleigh in the adapted basis are given by (see \cite{celledoni} for more details) \begin{equation}\label{EPSChaplygin}
\dot{v}^1=-\frac{a\sqrt{m}}{J+ma^2}v^1v^2,\quad\dot{v}^2=\frac{a\sqrt{m}}{J+ma^2}(v^1)^2, \end{equation}together with the admissibility conditions \begin{equation}\label{admisibilityCh}
\dot{x}_1=\frac{\cos\theta}{\sqrt{m}}v^2,\,\dot{x}_2=\frac{\sin\theta}{\sqrt{m}}v^{2},\,\dot{\theta}=\frac{1}{\sqrt{J+ma^2}}v^{1}
\end{equation}and the nonholonomic constraint $v^3=0$.

\subsection{Example: The nonholonomic particle}\label{example}

Consider a particle of unit mass evolving in $Q=\R^3$ with Lagrangian $\displaystyle{
L(x,y,z,\dot x,\dot y,\dot z)=\frac{1}{2}(\dot x^2+\dot y^2+\dot z^2)}$, and subject to the constraint $\dot x+y\,\dot z=0$. The nonholonomic system is defined by the annihilation of the one form $\mu(x,y,z)=(1,0,y)$. We denote $q(t)=(x(t),y(t),z(t))^T$ the vector of positions and $v(t)=(v_x(t),v_y(t),v_z(t))^T$ the corresponding vector of velocities.

 The distribution $\mathcal{D}$ is determined by  $\mathcal{D} = \hbox{span} \{Y_1, Y_2\}= \hbox{span}\Big{\{} \frac{\partial}{\partial y}, \frac{\partial}{\partial z}- y \frac{\partial}{\partial x}\Big{\}} $. Then, $\mathcal{D}^{\perp}=\{\frac{\partial}{\partial x}+y\frac{\partial}{\partial z}\}$. 
    
Let $(x,y,z,v^1,v^2)$ be induced coordinates on $\mathcal{D}$. Given the  vector fields $Y_1$ and $Y_2$ generating the distribution $\mathcal{D}$ we obtain the relations for $q\in \R^3$ given by $\displaystyle{Y_j(q)=\rho_j^1(q)\frac{\partial}{\partial x}+\rho_j^2(q)\frac{\partial}{\partial y}+\rho_j^3(q)\frac{\partial}{\partial z}}$, $j=1,2$. Then, $\rho_{1}^{1}= \rho_{1}^{3}=\rho_2^2=0, \quad\rho_{1}^{2}=\rho_2^3=1,\quad\rho_{2}^{1}=-y$.

 Each element $e\in {\mathcal D}_q$ is expressed as a linear combination of these vector fields: $e=v^1 Y_1(q)+ v^2 Y_2(q)$. Therefore, the vector subbundle $\tau_{\mathcal{D}}:{\mathcal D}\rightarrow \R^3$ is locally described by the coordinates $(x,y,\theta; v^1, v^2)$; the first three for the base and the last two, for the fibers.

Observe that 
$\displaystyle{
e=v^1\frac{\partial}{\partial y}+ v^2\left(\frac{\partial}{\partial z}-y\frac{\partial}{\partial x}
\right)}$ and, in consequence, ${\mathcal D}$ is described by the conditions (admissibility conditions): $
\dot{x}=-yv^2,\,\dot{y}=v^1,\,\dot{z}=v^2$ as  a vector subbundle of $TQ$ where $v^1$ and $v^2$ are the velocities relative to the basis of $\mathcal{D}$.

The nonholonomic bracket given
by $\lcf\cdot,\cdot\rcf=\mathcal{P}([\cdot,\cdot])$ satisfies $$\lcf
Y_1,Y_2\rcf=\mathcal{P}[Y_1,Y_2]=\mathcal{P}\left(-\frac{\partial}{\partial x}\right)=\frac{y}{1+y^2}\left(\frac{\partial}{\partial z}-y\frac{\partial}{\partial x}\right).$$

  Therefore, by using \eqref{relation} all the Christoffel symbols for the connection $\displaystyle{\nabla^{\mathcal{G}^{\mathcal{D}}}}$ vanish except $\Gamma_{12}^{2}$ which is given by $\displaystyle{\Gamma_{12}^{2}=\frac{y}{1+y^2}}$.

The restriction of the Lagrangian function $L$ on $\mathcal{D}$ in the adapted coordinates $(v^1,v^2)$ is given by $$\ell(x,y,z, y_1,y_2)=\frac{1}{2}\left((v^1)^2+(v^2)^2(y^2+1)\right).$$  Therefore, the nonholonomic equations for the constrained particle are given by
   \begin{equation}\label{eqq3}
     \dot{v}^1= 0,\quad
     \dot{v}^2= -\frac{y}{1+y^2}v^1v^2
   \end{equation} together with  the admissibility conditions $\dot{x}= -y v^2$,
     $\dot{y}= v^1$ and 
     $\dot{z}= v^2.$
      Then these equations define a time-continuous flow $F_t^{\mathcal{D}}:\mathcal{D}\Flder \mathcal{D}$, i.e. $F_t^{\mathcal{D}}((q(0),v(0)))=(q(t),v(t))$, where $q(t)=(x(t),y(t),z(t))^T$ and $v(t)=(v_1(t),v_2(t))^T$, $(q(0),v(0))\in\mathcal{D}$. 
      
      Note that only by taking an adapted basis of vector fields in the nonholonomic distribution $\mathcal{D}$), we reduced the quantity of equations to solve, without the needed to use Lagrange multipliers to enforce the nonholonomic constraint.

\section{Optimal trajectory tracking problem}\label{section3}

Next we present the tracking problem for nonholonomic systems as an optimal control problem. The objective is the tracking of a suitable reference trajectory $\gamma_r(t)$ for a mechanical system with velocity constraints as described in the previous section. It is assumed that $\gamma_r(t)\in\mathcal{D}$. 

We will analyze the case when the dimension of the inputs set, i.e., control distribution, is equal to the rank of $\mathcal{D}$. If the rank of $\mathcal{D}$ is equal to the dimension of the control distribution, the system will be called a \textit{fully actuated nonholonomic system}.

\begin{definition}
A \textit{solution of a fully actuated nonholonomic problem} is an
admissible curve $\gamma:I\rightarrow\mathcal{D}$ such that
$$\nabla_{\gamma(t)}^{\mathcal{G}^{\mathcal{D}}}\gamma(t)+grad_{\mathcal{G}^{\mathcal{D}}}V(\tau_{\mathcal{D}}(\gamma(t)))\in \Gamma (\tau_D),$$
or, equivalently, 
$$\nabla_{\gamma(t)}^{\mathcal{G}^{\mathcal{D}}}\gamma(t)+grad_{\mathcal{G}^{\mathcal{D}}}V(\tau_{\mathcal{D}}(\gamma(t)))=u^{A}(t)e_{A}(\tau_{\mathcal{D}}(\gamma(t))),$$
where $u^{A}$ are the control inputs.
\end{definition}

Locally, the above equations are given by 
\begin{eqnarray}
\dot{q}^{i}&=&\rho_{A}^{i}v^{A}\label{control1}\\
\dot{v}^{A}&=&-\Gamma_{CB}^{A}v^{C}v^{B}-(\mathcal{G}^{\mathcal{D}})^{AB}\rho_{B}^{i}(q)\frac{\partial V}{\partial q^{i}}+u^{A}.\label{control2}
\end{eqnarray}

As we mentioned in the Introduction, for trajectory tracking, the usual approach of stabilization of error dynamics \cite{kodi}, \cite{aradhana1}, \cite{aradhana2}, \cite{sanyal}  cannot be utilized for nonholonomic systems because the closed loop trajectory violates Brockett's condition. A common approach to trajectory tracking for nonholonomic systems found in the literature is the backstepping procedure \cite{Ha}, \cite{nijmeijer97}. This approach is done on basis of concrete examples, in particular, mobile robots or unicycle models.  In \cite{Ha}, \cite{nijmeijer97} the error dynamics of the unicycle model is shown to be in strict feedback form. Thereafter, integrator backstepping is employed to choose an appropriate Lyapunov function for stabilization of the error dynamics. This error dynamics does not evolve on the constrained manifold (unlike our approach). Therefore, Brockett's condition is not violated. However, since $\rho^i_A(q)$ is unknown in a general framework (i.e., they depend on the distribution determined in each particular case), the approach can not be generalized to solve the tracking problem for a general nonholonomic system with our method and then backstepping needs to be studied for each system.  So we propose a new approach by considering tracking problems for nonholonomic systems as optimal control problems, and we call this \textit{optimal trajectory tracking}.

In the following, we shall assume that all the control systems under consideration are controllable in the configuration space, that is, for any two
points $q_0$ and $q_f$ in the configuration space $Q$, there exists
an admissible control $u(t)$ defined on the \textcolor{blue}{control set}  
$\mathcal{U}\subseteq\R^{n}$ such that the system with initial condition $q_0$
reaches the point $q_f$ at time $T$, \textcolor{blue}{where $\mathcal{U}$ is unbounded} (see \cite{Bl} for
more details, \textcolor{blue}{Section $7.2$}). 

Given a cost function $
\mathcal{C}:\mathcal{D}\times \mathcal{U}\rightarrow\mathbb{R}$ the \textit{optimal control problem} consists of finding an admissible curve
$\gamma:I\rightarrow\mathcal{D}$ which is a solution of the fully actuated
nonholonomic problem given initial and final boundary conditions on
$\mathcal{D}$ and minimizing the cost functional
$$\mathcal{J}(\gamma(t),u(t)):=\int_{0}^{T}\mathcal{C}(\gamma(t),u(t))dt.$$

For trajectory tracking of a nonholonomic system we consider the following problem

\textbf{Problem (optimal trajectory tracking):} Given a reference trajectory $\gamma_{r}(t)=(q_r(t), v_r(t))$ on $\mathcal{D}$, find an admissible curve $\gamma(t)\in\mathcal{D}$, solving \eqref{control1}-\eqref{control2}, with prescribed boundary conditions on $\mathcal{D}$ and minimizing the cost functional 
\begin{align*}\mathcal{J}(\gamma(t)) =& \frac{1}{2}\int_{0}^{T} \left(||\gamma(t)- \gamma_{r}(t)||^2 + \epsilon ||u^A||^2 \right)\,dt+\omega\Phi(\gamma(T))\\
& =\frac{1}{2}\int_{0}^{T} \left(||q^{i}(t)- q^{i}_{r}(t)||^2 +||v^{A}(t)-v^{A}_{r}(t)||^2 +\epsilon ||u^A||^2 \right)\,dt+\omega\Phi(T, \gamma(T))
\end{align*} where $\epsilon>0$ is a regularization parameter, $\Phi:TQ\to\mathbb{R}$ is a terminal cost (Mayer term) and $\omega>0$ is a weight for the terminal cost. $\mathcal{C}$ and $\Phi$ are assumed to be continuously differentiable functions, and the final state $\gamma(T)$ is required to fulfill a constraint $r(\gamma(T), \gamma_{r}(T))=0$ with $r:\mathcal{D}\times\mathcal{D}\to\mathbb{R}^{d}$ and $\gamma_r\in\mathcal{D}$ given. The interval length $T$ may either be fixed, or appear as degree of freedom in the optimization problem. In this work we restrict ourselves to the case when $T$ is fixed.

\begin{remark}\label{rksingular}
Note that if $\epsilon=0$ then the optimal control problem turns into a singular optimal control problem (see \cite{OTT} Section $3.2$). This situation will be analyzed in a future work.
\end{remark}

\section{Conditions for optimality}\label{sec4}
In this section we derive necessary conditions for extrema in the optimal trajectory tracking problem. We present two approaches: the first one is based on the Hamiltonian point of view by considering Pontryagin's maximum principle, and  the second one is based on considering a Lagrangian point of view. In the Lagrangian approach, necessary conditions for extrema are derived as solutions of Euler-Lagrange equations for a Lagrangian defined as the cost functional for the optimal trajectory tracking problem. As we commented in the Introduction, the motivation to study the Lagrangian approach comes from the fact that by considering a Hamiltonian formalism, when we simulate the behavior of the planned trajectories by employing a classical integrator scheme in Section \ref{examplebad}, we can not obtain results that preserve the original qualitative structure of solutions. That is, despite we can reach the desired trajectory at the final time, the planned trajectories does not respect the original movements and behaviors of the continuous-time system,  and therefore  the construction of structure preserving numerical methods for this problem is needed. We construct such a structure preserving methods by discretizing the variational principle that we present in this section for the Lagrangian approach of the problem.

\subsection{Pontryagin Maximum Principle (PMP)}

In this section we apply Pontryagin's maximum principle to the optimal tracking problem. 

The Hamiltonian for the problem $\mathcal{H}:T^{*}\mathcal{D}\times \mathcal{U}\to\mathbb{R}$ is given by
\begin{align}\label{hamiltonian} \mathcal{H}(q,v, \lambda,\mu,u)=& \textcolor{blue}{\lambda_0\mathcal{C}(q^i,v^A,u^A)} + \lambda_i \rho_{A}^{i}(q)v^{A}+ \mu_A \dot{v}^A(q^{i},v^{A},u^{A}),
\end{align} where $\dot{v}^{A}$ comes from equation \eqref{control2} \textcolor{blue}{and 
$\lambda_0\geq 0$ is a fixed positive constant}. Note that $\lambda_i$ and $\mu_A$ are the costate variables. The second and third terms in \eqref{hamiltonian} correspond with the nonholonomic dynamics given in equations \eqref{eqq1} and \eqref{eqq2} paired with the costate variables.

We proceed as is usual in the literature (see for instance \cite{Bl} pp. $337$). \textcolor{blue}{We first restrict ourselves to the case of normal extremals, i.e., $\lambda_0\neq 0$}. The 
optimal curves $(q(t),v(t),\lambda(t),\mu(t), u^{\star}(t))$ must satisfy equations \eqref{control1} and \eqref{control2} together with the adjoint (or costate)  equations for $\mathcal{H}$, that is,  $$
-\dot{\lambda}_i=\frac{\partial \mathcal{H}}{\partial q^{i}}\hbox{ and }
-\dot{\mu}_A=\frac{\partial\mathcal{H}}{\partial v^{A}},
$$ where $u^{\star}$ satisfies, \begin{equation}\label{optimaleq}\mathcal{H}(q(t),v(t),\lambda(t),\mu(t),u^{\star}(t))=\min_{u\in\mathcal{U}}\mathcal{H}(q(t),v(t),\lambda(t),\mu(t), u).\end{equation} 

Given that $u^{\star}$ minimizes $\mathcal{H}$, then $u^{\star}$ is a critical point for $\mathcal{H}$ and may be determined by the condition \begin{equation}\label{stationary}\frac{\partial\mathcal{H}}{\partial u}(q(t),v(t),\lambda(t),\mu(t),u^{\star}(t))=0,\quad t\in[0,T].\end{equation} 

Note that by definition of $\mathcal{J}$,  $u^{\star}$ is determined
uniquely from the previous condition by the implicit function theorem. It follows that there exists a function
$\kappa$ such that $u^{\star}(t)=\kappa(q(t),v(t),\lambda(t),\mu(t))$. Then if $u^{\star}$ is defined implicitly as a function of $(q(t),v(t),\lambda(t),\mu(t))\in T^{*}\mathcal{D}$, by  equation \eqref{optimaleq} we can define the Hamiltonian function $\mathcal{H}^{*}:T^{*}\mathcal{D}\to\mathbb{R}$ by $$\mathcal{H}^{*}(q(t),v(t),\lambda(t),\mu(t))=\mathcal{H}(q(t),v(t),\lambda(t),\mu(t), u^{\star}(t)).$$ $\mathcal{H}^{*}$ defines a Hamiltonian vector field $X_{\mathcal{H}^{*}}$ on $T^{*}\mathcal{D}$ with respect to the canonical symplectic structure on $T^{*}\mathcal{D}$ given by $\omega_{\mathcal{D}}=dq^i\wedge d\lambda_i+dv^A\wedge d\mu_A$.

The PMP applied to our particular problem, together with the constraints induced by the terminal cost and the boundary conditions gives the following necessary conditions:

\begin{itemize}
\item[(i)] Stationary condition: from \eqref{stationary}, $\mu_{A}=-\textcolor{blue}{\lambda_0}\epsilon u^{A},$ that is, $\displaystyle{(u^A)^{\star}=-\frac{\mu_A}{\textcolor{blue}{\lambda_0}\epsilon}}$.
\item[(ii)] State equations: Equations \eqref{control1} and \eqref{control2}, with $u^A$ determined by the stationary condition.
\item[(iii)] Adjoint equations (or costate equations): \begin{align*}
-\dot{\lambda}_i=&\frac{\partial \mathcal{H}^{*}}{\partial q^{i}}=\textcolor{blue}{\lambda_0}(q^{i}-q^{i}_{r})+\lambda_j\frac{\partial\rho_{A}^{j}(q)}{\partial q^{i}}v^{A}+\mu_A\frac{\partial\dot{v}^{A}}{\partial q^{i}},\\
-\dot{\mu}_A=&\frac{\partial\mathcal{H}^{*}}{\partial v^{A}}=\textcolor{blue}{\lambda_0}(v^{A}-v_r^{A})+\lambda_i\rho_A^{i}(q)+\mu_B\frac{\partial \dot{v}^{B}}{\partial v^{A}}\,
\end{align*}
\item[(iv)] Constraint induced by terminal \textcolor{blue}{condition}: $r(\gamma(T), \gamma_{r}(T))=0$,
\item[(v)] \textcolor{blue}{Transversality conditions: $\gamma(0):=(q(0),v(0))\in\mathcal{D}$, 
 \begin{eqnarray*}
 \lambda_i(T)&=&\omega\frac{\partial\Phi}{\partial q^{i}}(T, \gamma(T))+\lambda_T \frac{\partial r}{\partial q^i}(\gamma(T), \gamma_r(T)),\\
\mu_A(T)&=&\omega\frac{\partial\Phi}{\partial v^{A}}(T, \gamma(T))+\lambda_T \frac{\partial r}{\partial v^A}(\gamma(T), \gamma_r(T)).
\end{eqnarray*}
}
\end{itemize}


\textcolor{blue}{
Observe that the solutions of the optimal control problem are the critical points of the functional
\begin{eqnarray*}
{\widetilde{J}}(\gamma, u, \lambda, \mu,  \lambda_T)&=&\omega \Phi(T, \gamma(T))+\lambda_T r(\gamma(T), \gamma_r(T))\\
&& +\int_0^T\left[ 
\lambda_0\mathcal{C}(q(t),v(t),u(t)) + \lambda_i(t) (\dot{q}^i(t)- \rho_{A}^{i}(q(t))v^{A}(t))\right.\\
&&\qquad \left.+ \mu_A(t) (\dot{v}^A(t)-\dot{v}^A(q(t),v(t),u(t)))\right]\,dt,
\end{eqnarray*}
with $\omega>0, \lambda_0\geq 0$, $\gamma(0)\in {\mathcal D}$, $\lambda_T\in\mathbb{R}$ and $\gamma_r:[0, T]\rightarrow {\mathcal D}$ given. 
}

\textcolor{blue}{Note that in the abnormal case, that is, when $\lambda_0=0$, it follows that $\mu_A=0$ and therefore the adjoint equations become in}
 \begin{equation*}
\textcolor{blue}{-\dot{\lambda}_i=\lambda_j\frac{\partial\rho_{A}^{j}(q)}{\partial q^{i}}v^{A} ,\quad
0=\lambda_i\rho_A^{i}(q).}
\end{equation*}

\begin{remark}
\textcolor{blue}{In the situation for the study of abnormal solutions, the necessary conditions cannot use the information of the cost function $\mathcal{C}$ to select minimizers. That is, abnormal solutions are not useful solutions for our trajectory tracking problem, since the problem formulation for optimal trajectory tracking depends explicitly in the distance between the desired trajectory and the optimal one. The unique  condition that we need in our work is the existence of normal solutions, which in our case,  are guaranteed by assuming the controllability of the linearized state equations (see  \cite{mala}). This is the typical controllability hypothesis  assumed for  trajectory tracking and it is the general case in control nonholonomic dynamics.}
\end{remark}

\subsection{Example: Optimal trajectory tracking for the Chaplygin sleigh}\label{chch2}
Consider the Chaplyigin sleigh of Example \ref{chaplygin} but subject to input controls. These control inputs are denoted by $u_1$ and $u_2.$ The first control input
corresponds to a force applied perpendicular to the center of mass
of the sleigh and the second control input corresponds to the torque applied about the
vertical axis.

The controlled Euler-Lagrange
equations are given by
\begin{equation}\label{eqchcontrol2}
\dot{v}^1=-\frac{a\sqrt{m}}{J+ma^2}v^1v^2+u_1,\quad\dot{v}^2=\frac{a\sqrt{m}}{J+ma^2}(v^1)^2+u_2.
\end{equation}
together with the admissibility conditions \begin{equation}\label{eqadm}
\dot{x}_1=\frac{\cos\theta}{\sqrt{m}}v^2,\,\dot{x}_2=\frac{\sin\theta}{\sqrt{m}}v^{2},\,\dot{\theta}=\frac{1}{\sqrt{J+ma^2}}v^{1}
\end{equation}and the nonholonomic constraint $v^3=0$.

Let $\gamma_{r}(t)=((x_1)_{r}(t),(x_2)_{r}(t),\theta_r(t),v^{1}_r(t),v_{r}^{2}(t))$  be the reference trajectory, which follows the constraint $v^{3}_r=0$ for all time $t$ and the dynamical equations for the Chaplygin sleigh. \textcolor{blue}{ In this case, we assume that the final cost is  $\Phi(T, \gamma(T))=0$, and the  constraint $r(\gamma(T),\gamma_r)$ is given by \begin{align*}
r(\gamma(T),\gamma_r)=&|x_1(T)-(x_1)_r(T) |^2 + |x_2(T)- (x_2)_r(T) |^2 +|\theta(T)- \theta_r(T) |^2 \\&+ |v^1(T)- v^{1}_r(T)|^2 + |v^2(T)- v^{2}_{r}(T)|^2.\end{align*}     
  }
  
The Hamiltonian for the PMP is given by
\begin{align*}
\mathcal{H}(q,v, \lambda, \mu,u)= &\frac{\textcolor{blue}{\lambda_0}\epsilon}{2}(u_1^2+u_2^2)+\frac{\textcolor{blue}{\lambda_0}}{2}(|x_1- (x_1)_r |^2 + |x_2- (x_2)_r |^2 +|\theta- \theta_r |^2 \\&+ |v^1- v^{1}_r|^2 + |v^2- v^{2}_{r}|^2 )+\lambda_1\frac {\cos\theta}{\sqrt{m}}v^2+\lambda_2\frac{\sin\theta}{\sqrt{m}}v^2+\frac{\lambda_3v^1}{\sqrt{J+ma^2}}\\&+\mu_1\left(u_1-\frac{a\sqrt{m}}{J+ma^2}v^1v^2\right)+\mu_2\left(u_2+\frac{a\sqrt{m}}{J+ma^2}(v^1)^2\right)\\
\end{align*}

Note that, $\displaystyle{u_1^{\star} = -\frac{\mu_1}{\textcolor{blue}{\lambda_0}\epsilon}\; \hbox{ and } \; u_2^{\star} = -\frac{\mu_2}{\textcolor{blue}{\lambda_0}\epsilon}}$. Therefore denoting by $q=(x_1,x_2,\theta)$, the optimal Hamiltonian $\mathcal{H}^{*}$ is given by 
\begin{align*}
\mathcal{H}^{*}(q,v, \lambda,\mu)= &\frac{\textcolor{blue}{\lambda_0}}{2} \big\{|x_1- (x_1)_r |^2 + |x_2- (x_2)_r |^2 +|\theta- \theta_r |^2+ |v^1- v^{1}_r|^2 + |v^2- v^{2}_{r}|^2 \big \}\\&+\lambda_1\frac {\cos\theta}{\sqrt{m}}v^2+\lambda_2\frac{\sin\theta}{\sqrt{m}}v^2+\frac{\lambda_3v^1}{\sqrt{J+ma^2}}-\frac{\mu_1^{2}}{2\textcolor{blue}{\lambda_0}\epsilon}-\frac{\mu_2^2}{2\textcolor{blue}{\lambda_0}\epsilon}\\&-\mu_1\frac{a\sqrt{m}}{J+ma^2}v^1v^2+\mu_2\frac{a\sqrt{m}}{J+ma^2}(v^1)^2.
\end{align*}

The adjoint equations are  
\begin{align}\label{costateeqn}
\dot{\lambda}_1&= -\textcolor{blue}{\lambda_0}(x_1-(x_1)_r),\, \dot{\lambda}_2 = -\textcolor{blue}{\lambda_0}(x_2-(x_2)_r),\nonumber\\ 
\dot{\lambda}_3 &= \textcolor{blue}{\lambda_0}(\theta_r-\theta)+\lambda_1\frac{\sin\theta}{\sqrt{m}}v^2-\lambda_2\frac{\cos\theta}{\sqrt{m}}v^2,\\ \nonumber
\dot{\mu}_1 &= - \textcolor{blue}{\lambda_0}(v^1- v^{1}_{r})-\lambda_3\frac{1}{\sqrt{J+ma^2}}+\mu_1v^2\frac{a\sqrt{m}}{J+ma^2}-\mu_2v^1\frac{2a\sqrt{m}}{J+ma^2}, \nonumber\\
\dot{\mu}_2 &= -\textcolor{blue}{\lambda_0} (v^2- v^{2}_{r})-\lambda_1\frac{\cos\theta}{\sqrt{m}}-\lambda_2\frac{\sin\theta}{\sqrt{m}} +\mu_1v^1\frac{a\sqrt{m}}{J+ma^2}.\nonumber
\end{align} Finally, the state equations are given now by \begin{equation}\label{eqchcontrol}
\dot{v}^1=-\frac{a\sqrt{m}}{J+ma^2}v^1v^2-\frac{\mu_1}{\textcolor{blue}{\lambda_0}\epsilon},\quad\dot{v}^2=\frac{a\sqrt{m}}{J+ma^2}(v^1)^2-\frac{\mu_2}{\textcolor{blue}{\lambda_0}\epsilon}. \end{equation} together with the admissibility conditions \begin{equation}\label{eqadm2}
\dot{x}_1=\frac{\cos\theta}{\sqrt{m}}v^2,\,\dot{x}_2=\frac{\sin\theta}{\sqrt{m}}v^{2},\,\dot{\theta}=\frac{1}{\sqrt{J+ma^2}}v^{1}
\end{equation} 
\textcolor{blue}{In addition, the boundary  conditions and transversality conditions must be satisfied, in particular,  the optimal trajectory verifies that  $\gamma(T)$ matches exactly with $\gamma_r(T)$ .}%

\subsection{Example: Optimal trajectory tracking for the nonholonomic particle}\label{examplebad}
Consider the situation of Example \ref{example}. Let $\gamma_{r}=(x_r(t),y_r(t),z_r(t),v^{1}_{r},v^{2}_r)$  be the reference trajectory, which follows the constraint $\dot{x}_r= y_r \dot{z}_r$ for all time $t$ and the dynamical equations for the nonholonomic particle. We wish to control the velocity of the nonholonomic particle. To do that, we add control inputs in the fiber coordinates $v^{1}$ and $v^{2}$. Therefore the dynamical control system to study is given by \begin{equation}\label{eqqcontrol}
     \dot{v}^1= u^{1},\quad
     \dot{v}^2=u^{2} -\frac{y}{1+y^2}v^1v^2,
   \end{equation} together with  the admissibility conditions $\dot{x}= -y v^2$,
     $\dot{y}= v^1$ and 
     $\dot{z}= v^2.$
     

The Hamiltonian for the PMP is given by
\begin{align*}
\mathcal{H}(q,v, \lambda, \mu,u)= &\frac{\textcolor{blue}{\lambda_0}}{2} \left(|x- x_r |^2 + |y- y_r |^2 +|z- z_r |^2 +|v^1- v^{1}_r|^2 + |v^2- v^{2}_{r}|^2 \right.\\
+&\left. \epsilon (u^1)^2+\epsilon(u^2)^2 \right)-\lambda_1 y v^2 + \lambda_2 v^1 + \lambda_3 v^2 + \mu_1 u^1 \\&+ \mu_2 \left(u^2-\frac{y}{1+y^2}v^1v^2\right).
\end{align*}

Note that, $\displaystyle{u_1^{\star} = -\frac{\mu_1}{\textcolor{blue}{\lambda_0}\epsilon}\; \hbox{ and } \; u_2^{\star} = -\frac{\mu_2}{\textcolor{blue}{\lambda_0}\epsilon}}$. Therefore the optimal Hamiltonian $\mathcal{H}^{*}$ is given by 
\begin{align*}
\mathcal{H}^{*}(q,v, \lambda,\mu)= &\frac{\textcolor{blue}{\lambda_0}}{2} \big\{|x- x_r |^2 + |y- y_r |^2 +|z- z_r |^2+ |v^1- v^{1}_r|^2 + |v^2- v^{2}_{r}|^2 \big \}\\&-\lambda_1 y v^2 + \lambda_2 v^1 + \lambda_3 v^2-\frac{1}{2\textcolor{blue}{\lambda_0}\epsilon}(\mu_1^2+\mu_2^2)-\mu_2v^1v^2\frac{y}{1+y^2}.
\end{align*}

The adjoint equations are  
\begin{align}\label{costateeqn2}
\dot{\lambda}_1&= -\textcolor{blue}{\lambda_0}(x-x_r),\quad \dot{\lambda}_3 = -\textcolor{blue}{\lambda_0}(z-z_r),\nonumber\\ 
\dot{\lambda}_2 &= \lambda_1 v^2 -\textcolor{blue}{\lambda_0}(y-y_r)+v^1v^2\mu_2\left(\frac{y^2-1}{(y^2+1)^{2}}\right),\\ \nonumber
\dot{\mu}_1 &= -\lambda_2 - \textcolor{blue}{\lambda_0}(v^1- v^{1}_{r})-\mu_2\frac{y}{1+y^2} v^2,\\  \nonumber
\dot{\mu}_2 &= -\lambda_3+ \lambda_1 y - \textcolor{blue}{\lambda_0}(v_2- v_{2}^{r}) -\mu_2\frac{y}{1+y^2} v_1.
\end{align}  Finally, the state equations are given now by \begin{equation}\label{eqqcontrol22}
     \dot{v}^1= -\frac{\mu_{1}}{\textcolor{blue}{\lambda_0}\epsilon},\quad
     \dot{v}^2= -\frac{\mu_{2}}{\textcolor{blue}{\lambda_0}\epsilon} -\frac{y}{1+y^2}v^1v^2,
   \end{equation} together with  the admissibility conditions $\dot{x}= -y v^2$,
     $\dot{y}= v^1$ and 
     $\dot{z}= v^2$. 
 \textcolor{blue}{In addition, we consider a final cost $\Phi(T,\gamma(T))$ (but not a function $r$) and the boundary  conditions and transversality conditions must be satisfied.}

We now test with numerical simulations how the proposed method works. 
We choose an arbitrary trajectory satisfying the nonholonomic dynamics and we solve the boundary value problem by using a single shooting method.

Denote by $F_{\mu}^{\lambda}:[0,T]\times T^{*}\mathcal{D}\to T^{*}\mathcal{D}$ the integral flow given by equations \eqref{costateeqn2} on $T^{*}\mathcal{D}$  and $\gamma(0)\in\mathcal{D}$ the initial condition for the state dynamics. The initial guess for the initial condition of the costate variables  is denoted by $ \alpha=F_{\mu}^{\lambda}(0)$. We wish to find the initial condition of the costates for which $F_{\mu}^{\lambda}(T, \gamma(0),  \alpha)=(0_{1\times 5})^T$ . The goal is to find the root of the polynomial \begin{small}\[
F_{\mu}^{\lambda}(\alpha)= \begin{pmatrix}
{\lambda}_1(T, \gamma(0),\alpha)+\omega( x(T, \alpha)- x_r(T)) \\
{\lambda}_2(T, \gamma(0), \alpha) +\omega(  y(T, \alpha)- y_r(T))\\
{\lambda}_3(T, \gamma(0), \alpha) +\omega(  z(T, \alpha)- z_r(T))\\
{\mu}_1(T, \alpha)\\
{\mu}_2(T, \alpha)
\end{pmatrix}
\] \end{small}where $T\in\mathbb{R}^{+}$ is the final time, $\omega\in\mathbb{R}^{+}$ is a weight for the terminal cost and $F_{\mu}^{\lambda}(\tau, \gamma(0), p_0)$ is the flow of the adjoint equations \eqref{costateeqn2} starting at $(\gamma(0),p_0)$. The root finder used in both situations was the \textit{fsolve} routine in MATLAB. 

\textit{Case 1: Singular case}.

\begin{figure}[tbhp]\label{fig1}
    \includegraphics[width=7cm]{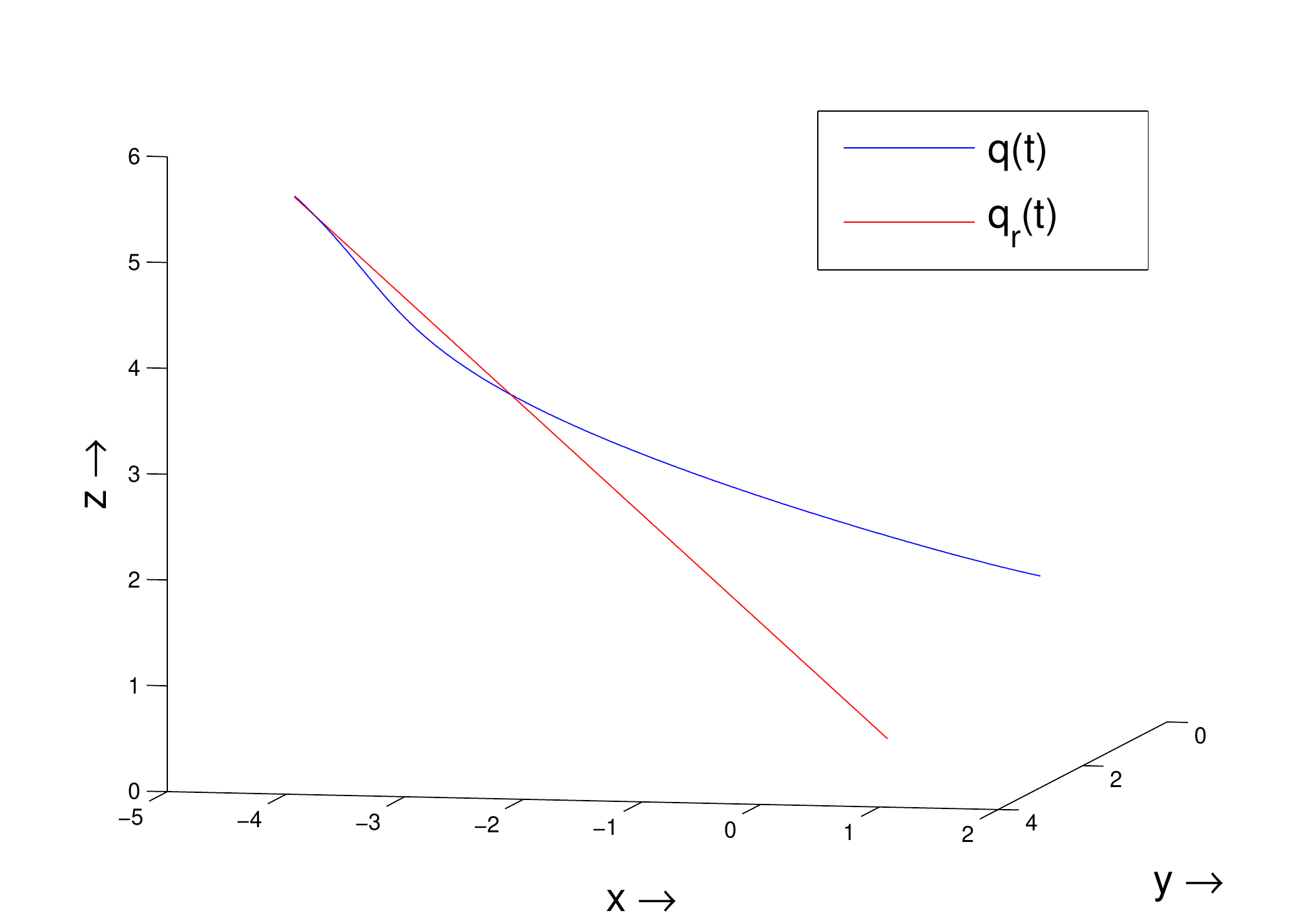}
    \hspace{-0.25cm}\includegraphics[width=7cm]{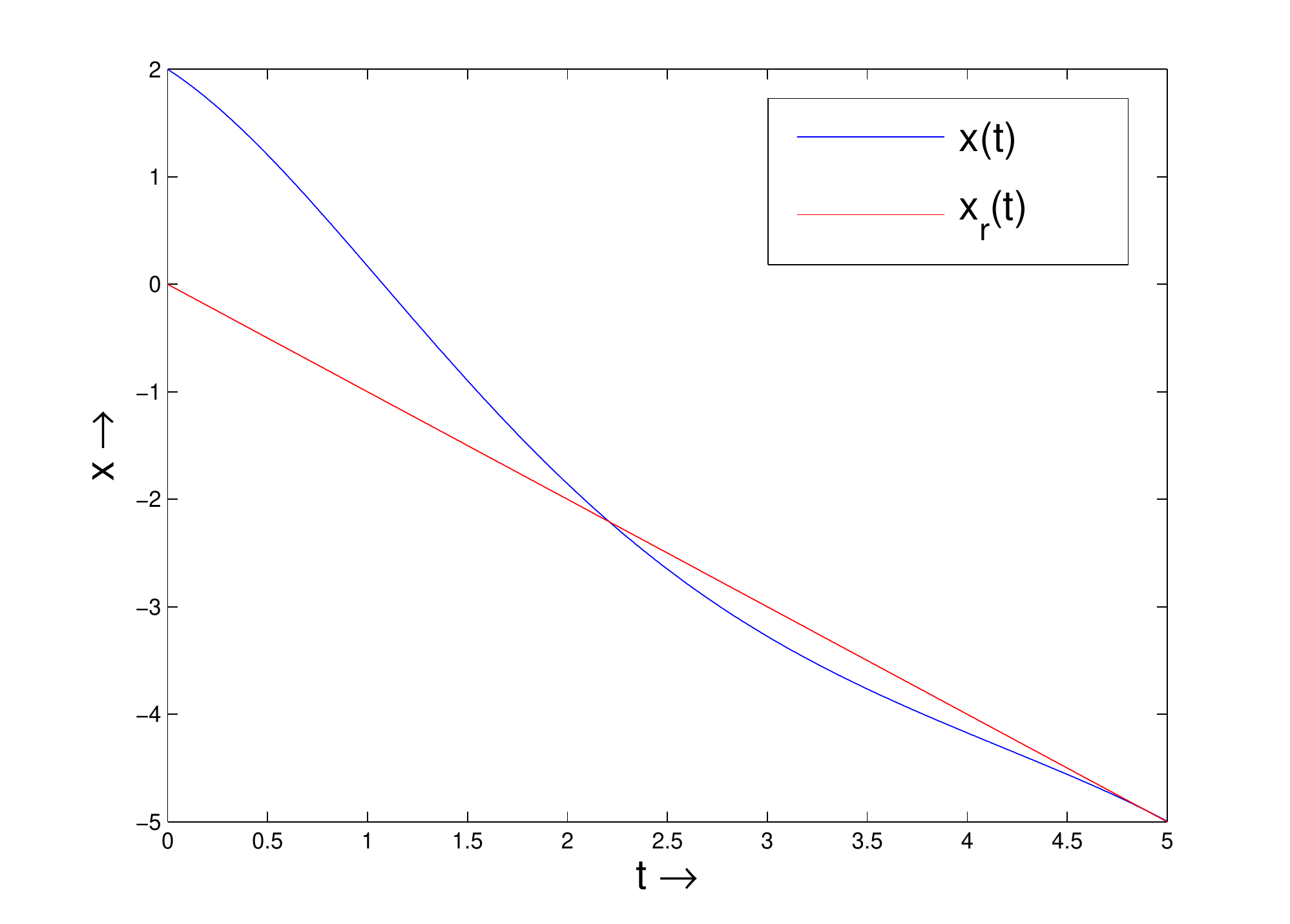}
    \newline
    \includegraphics[width=7cm]{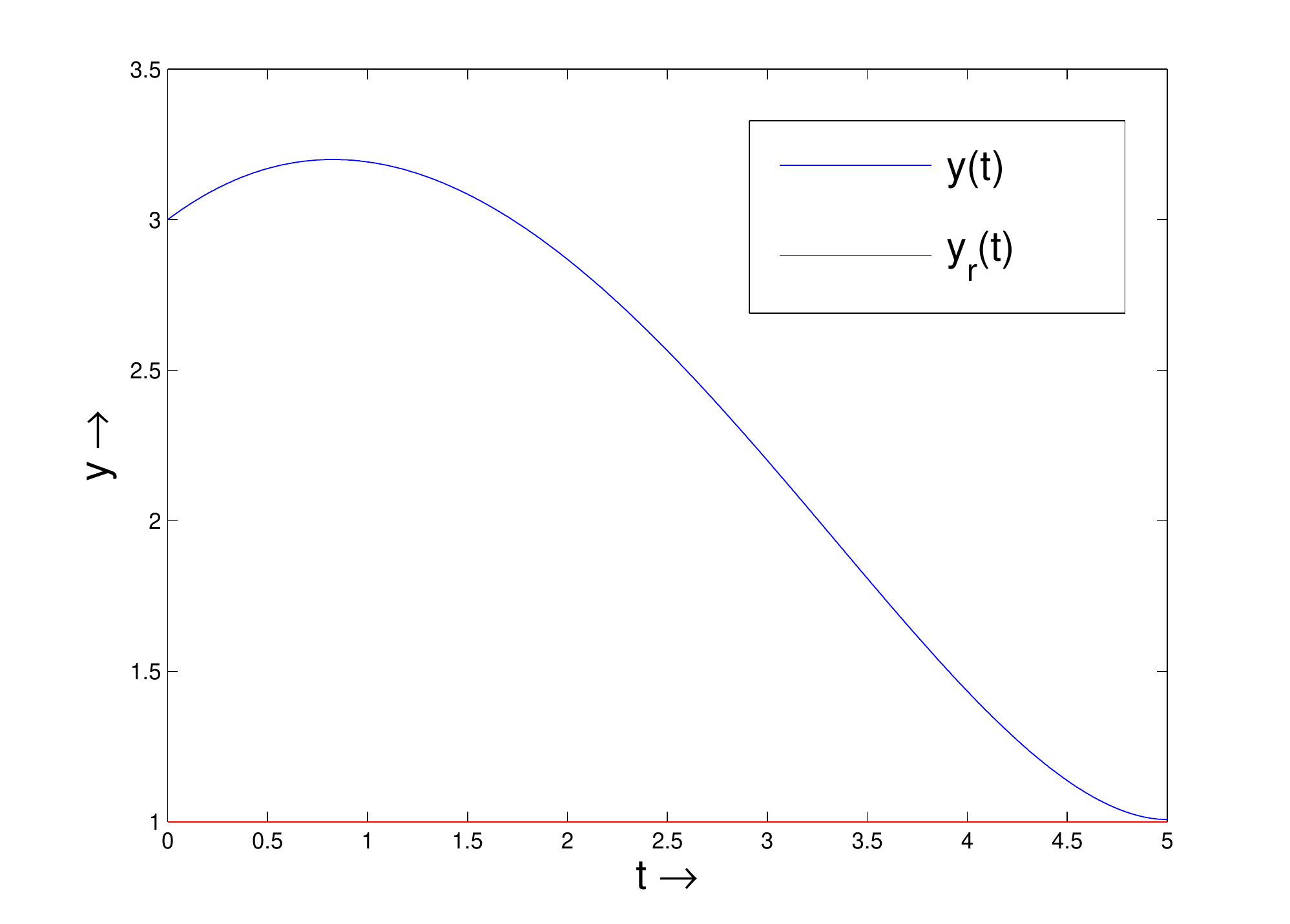}
    \hspace{-0.25cm}\includegraphics[width=7cm]{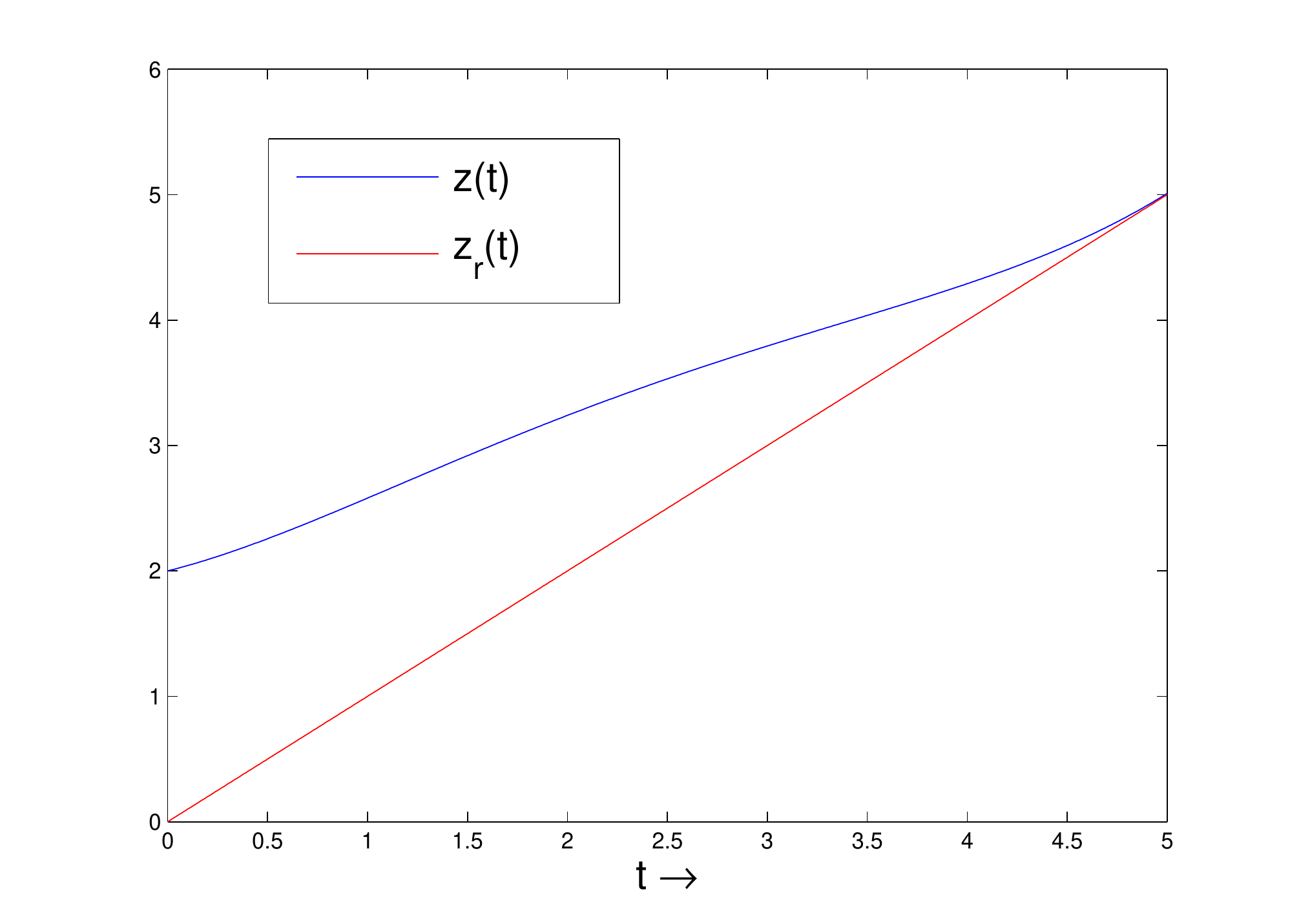}
    \newline
    \includegraphics[width=7cm]{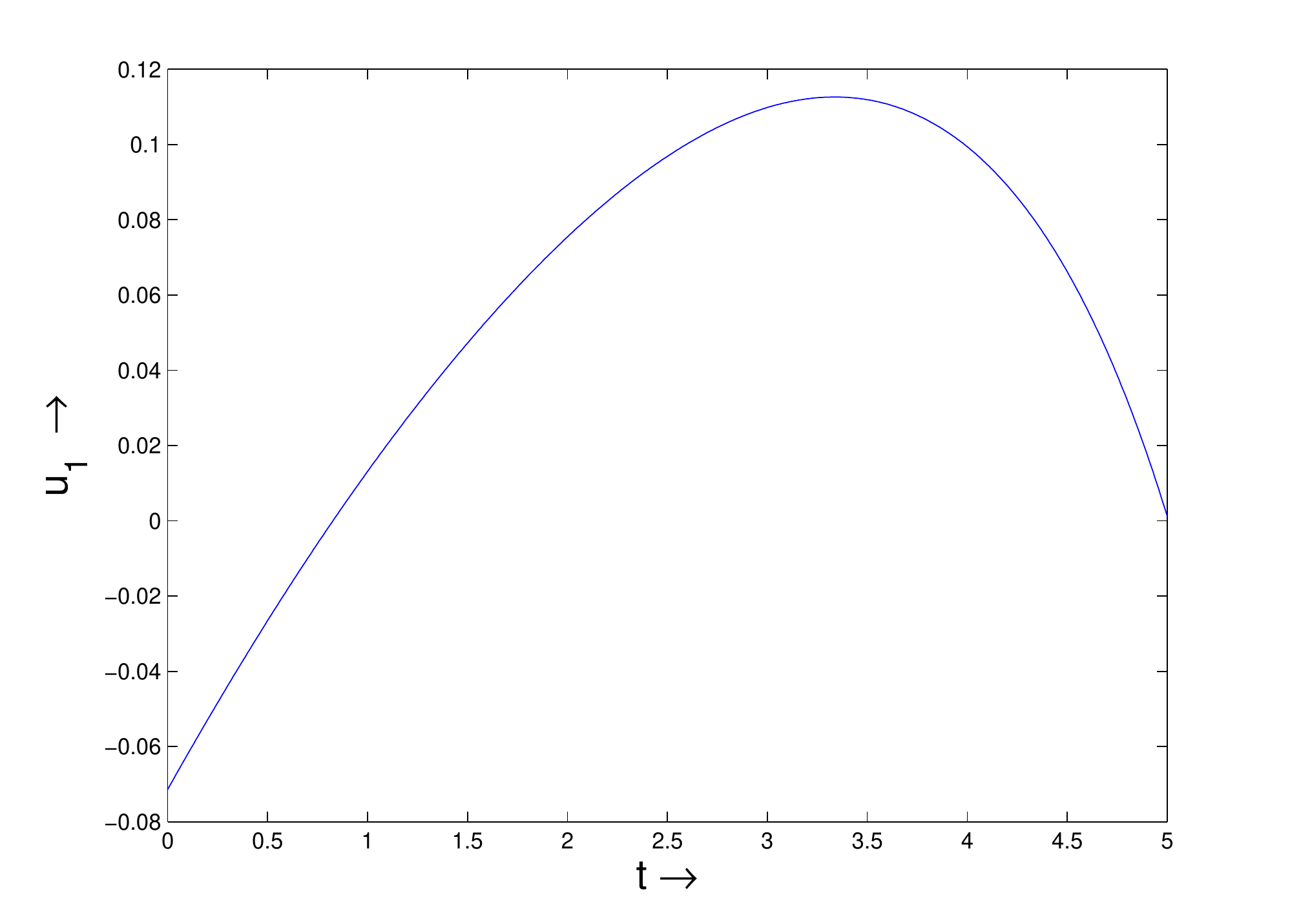}
    \hspace{-0.25cm}\includegraphics[width=7cm]{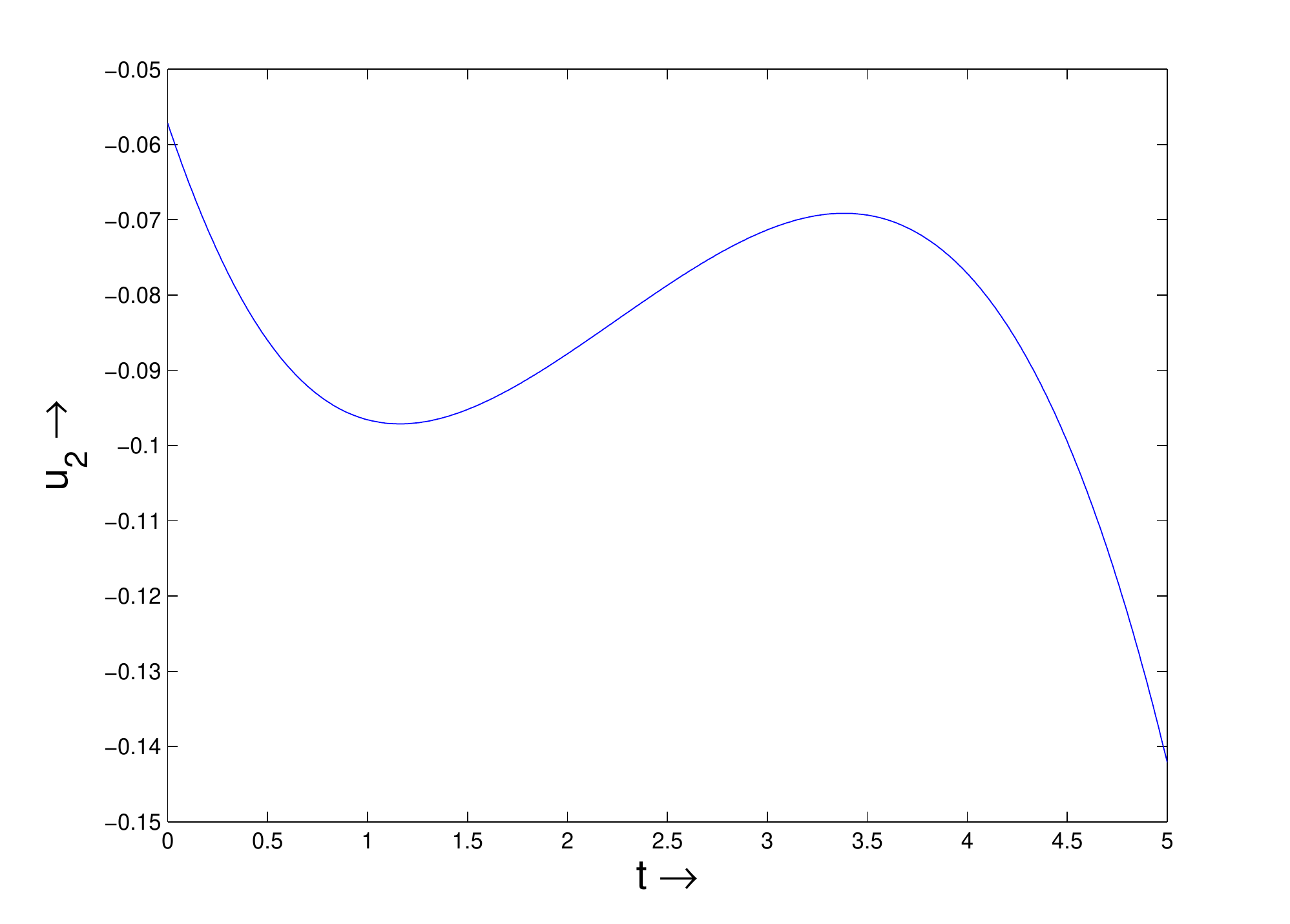}
  \caption{Singular case, $c_1=0$: Trajectories minimizing the cost function $\mathcal{J}$, evolving on $\mathcal{D}$ and tracking the reference trajectory $\gamma_r$ in time $T$ and control inputs}
  \end{figure}

For the initial condition $\gamma(0)= \begin{pmatrix}
2 & 3 & 2; & 0.5 & 0.4
\end{pmatrix}$ and reference trajectory $\gamma_r(t)= \begin{pmatrix}
 -t & 1 & t;& 0& 1
\end{pmatrix}$, $p_0= 0_{1\times 5}$, $T=5$, $\omega=1$ and $\epsilon=9$ we exhibit the results in Figure \ref{fig1}.

\textit{Case 2: Arbitrary reference trajectory}

\begin{figure}[h!]\label{fig3}
    \includegraphics[width=7cm]{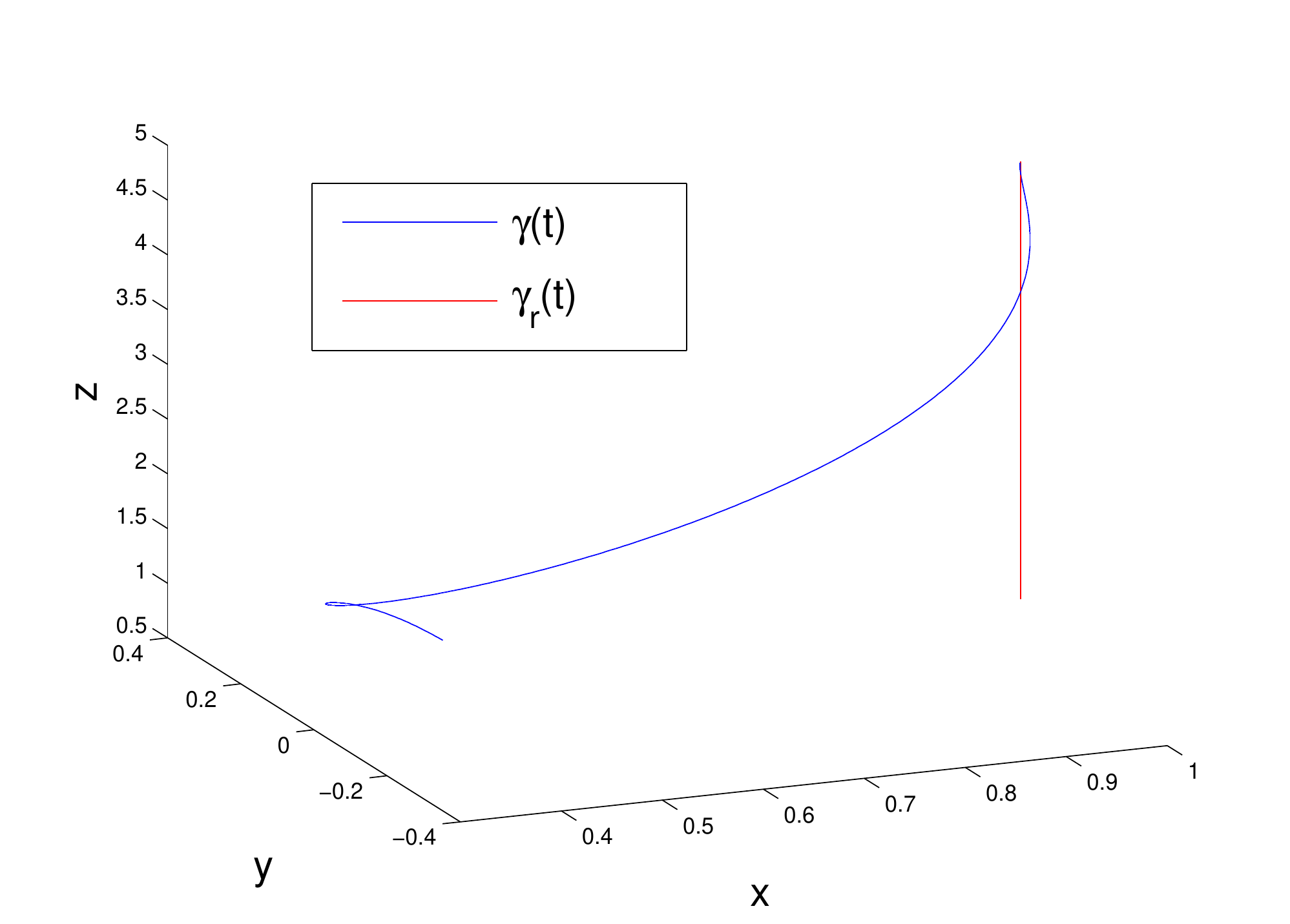}
    \hspace{-0.25cm}\includegraphics[width=7cm]{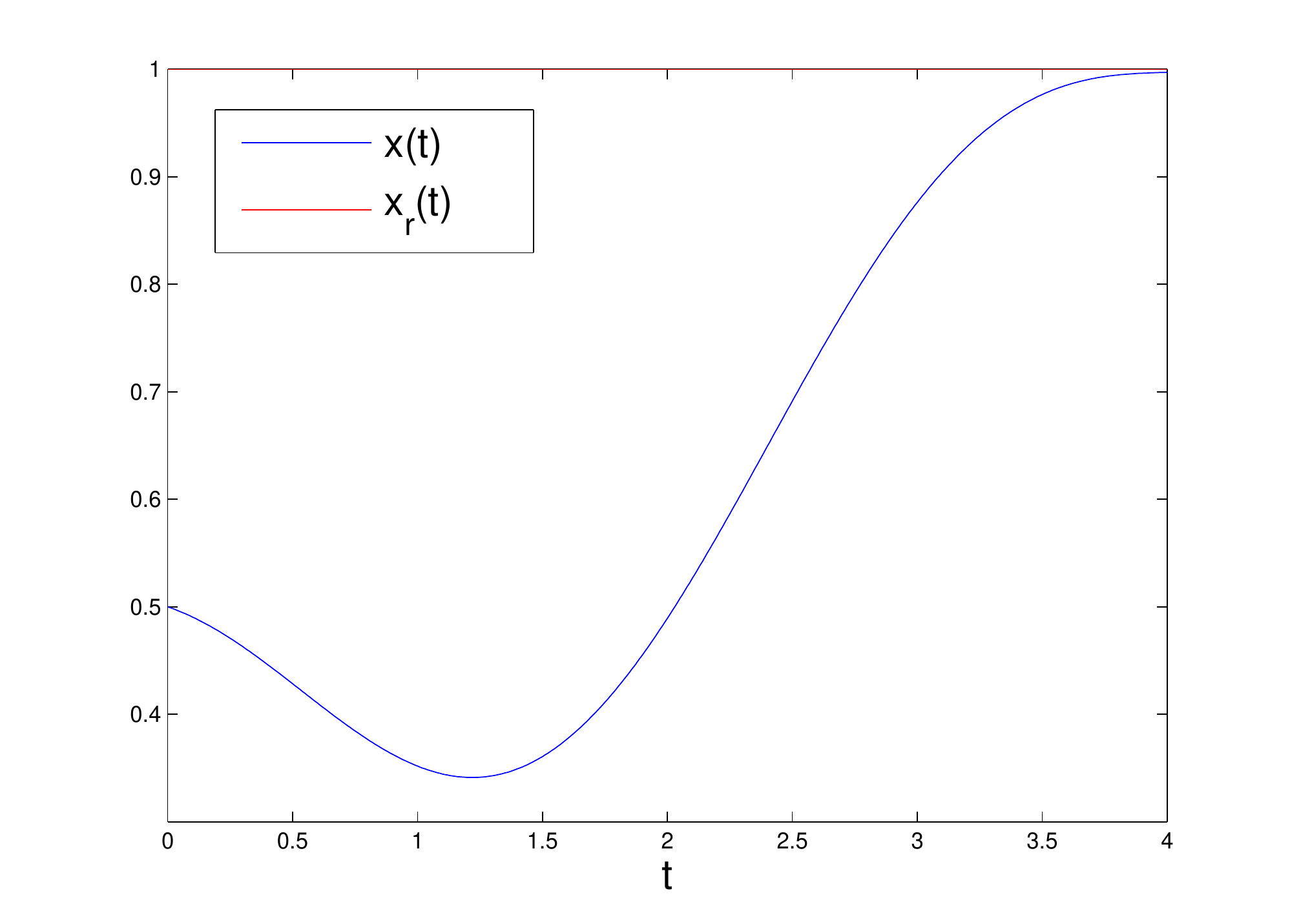}
    \newline
    \includegraphics[width=7cm]{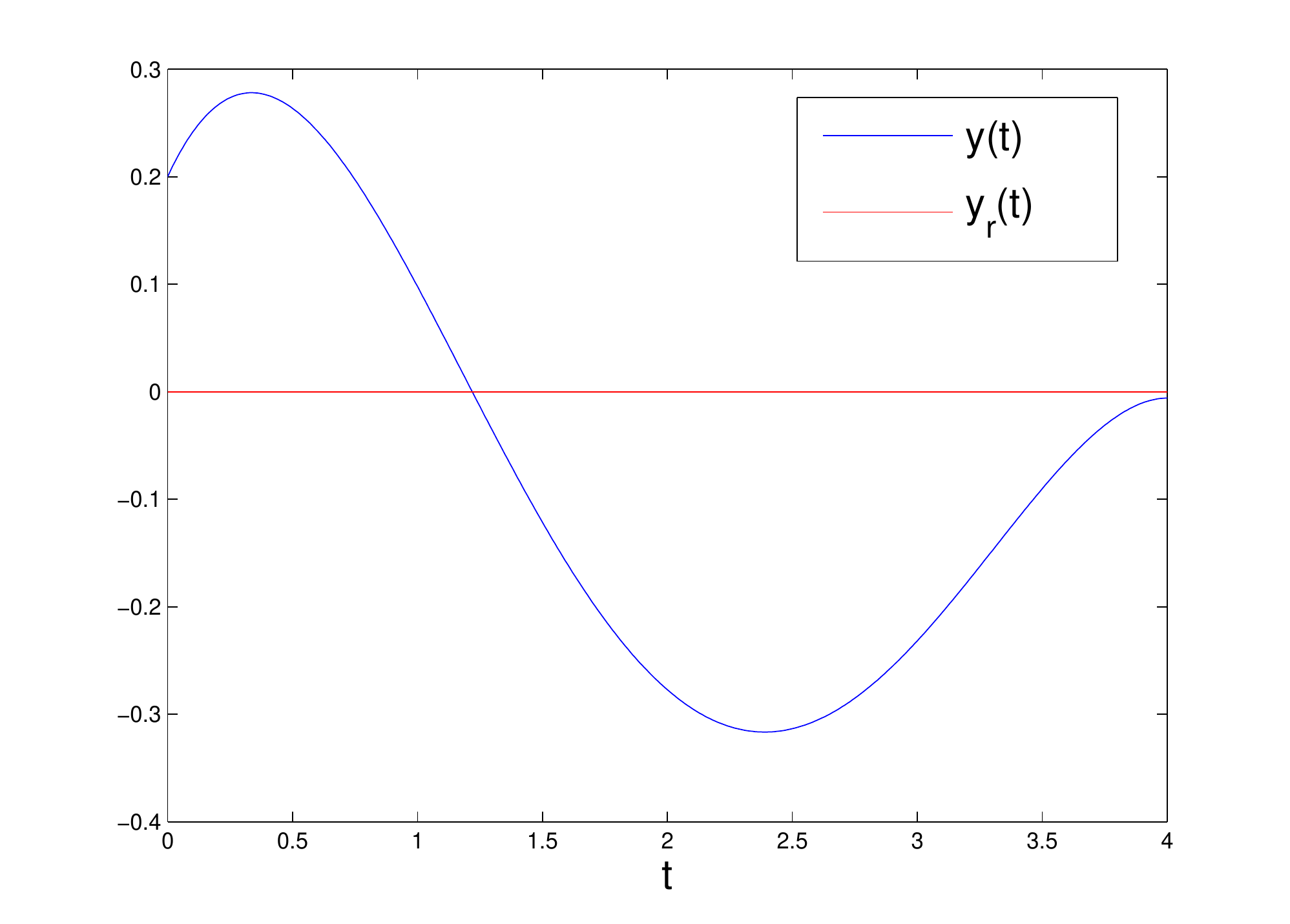}
    \hspace{-0.25cm}\includegraphics[width=7cm]{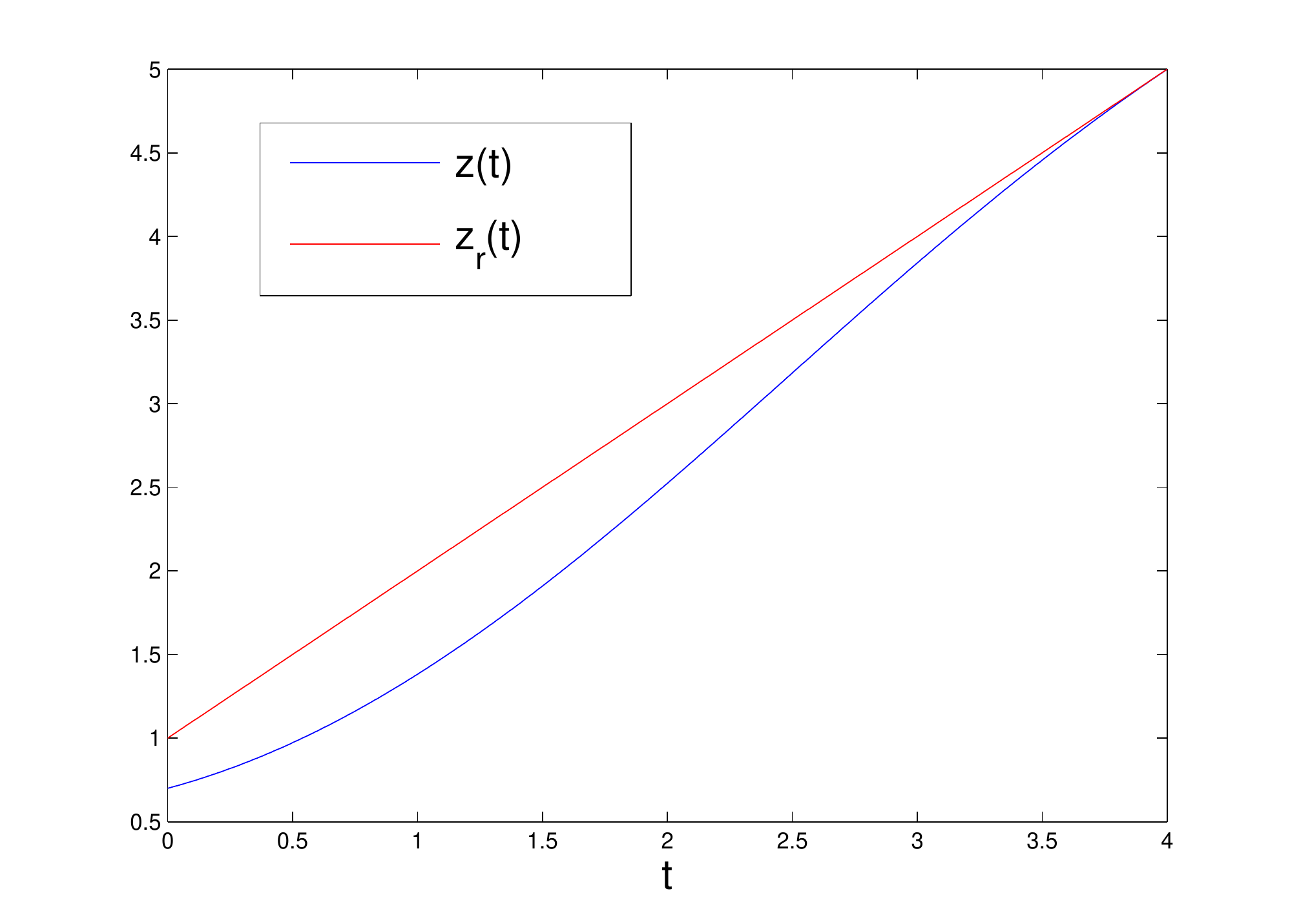}
    \newline
    \includegraphics[width=7cm]{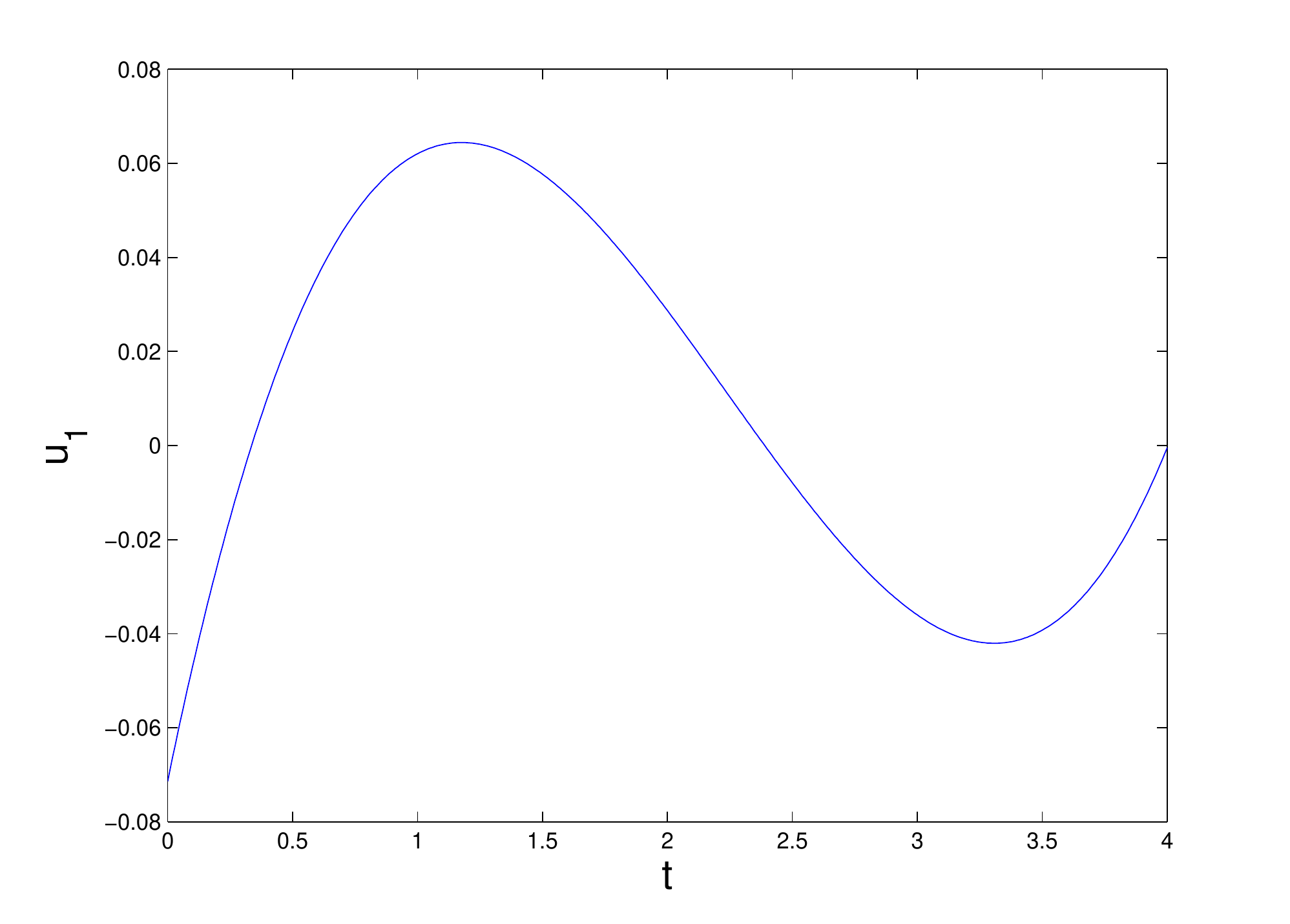}
    \hspace{-0.25cm}\includegraphics[width=7cm]{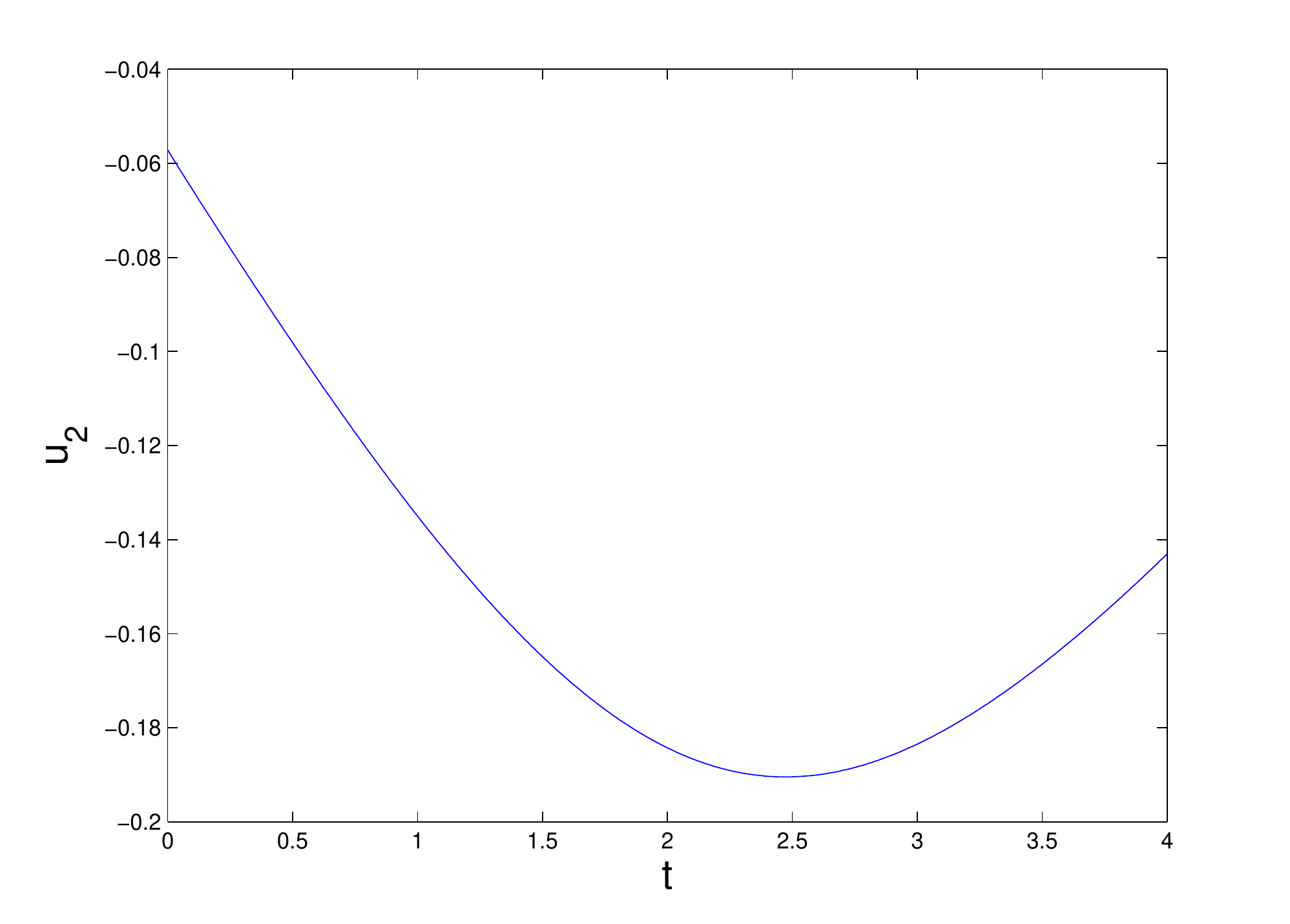}
  \caption{Trajectories minimizing the cost function $\mathcal{J}$, evolving on $\mathcal{D}$ and tracking the reference trajectory $\gamma_r$ in time $T$ and control inputs}
\end{figure}

For the intial condition $\gamma(0)=$ $\begin{pmatrix}
0.5 & 0.2 & 0.7; & 0.5 & 0.4
\end{pmatrix}$ and reference trajectory $\gamma_r(t)=$ $\left(
1, 0,  t+1,0,1\right) $, $p_0= 0_{1\times 5}$, $T=4$, $\omega=1$ and $\epsilon=7$ we exhibit the results in Figure \ref{fig3}.

Minimizing the cost functional, while evolving on the constraint submanifold and remaining differentiable by solving a boundary value problem using a single shooting method is a difficult task and not always numerically stable. Moreover, here we are not considering time as an independent variable, which will only complicate things further.  The need for using proper regularization parameters and final weights is crucial in order to get accurate results. In the next section we will improve the behavior in simulations  by constructing variational integrators. 

\subsection{Variational (Lagrangian) approach}
Next we derive necessary conditions for optimality in the optimal control problem following a variational approach as in \cite{BlCoGuMdD}, \cite{leo}, \cite{leo2}. Define the submanifold $\mathcal{D}^{(2)}$ of $T\mathcal{D}$ by $\mathcal{D}^{(2)}:=\{a\in T\mathcal{D}\mid a=\dot{\gamma}\}$, where $\gamma:I\rightarrow\mathcal{D}$ is an admissible curve.
 We can choose coordinates $(x^{i},v^{A},\dot{v}^{A})$ on
$\mathcal{D}^{(2)}$, where the inclusion on $T\mathcal{D}$,
$i_{\mathcal{D}^{(2)}}:\mathcal{D}^{(2)}\hookrightarrow
T\mathcal{D}$, is given by
$i_{\mathcal{D}^{(2)}}(q^{i},v^{A},\dot{v}^{A})=(q^{i},v^{A},\rho_{A}^{i}(q)v^{A},\dot{v}^{A})$. Therefore, $\mathcal{D}^{(2)}$ is locally described by the constraint on $T\mathcal{D}$ given by $\dot{q}^{i}-\rho_{A}^{i}v^{A}=0$.

The optimal control problem can be alternatively
studied by the function
$\mathcal{L}:\mathcal{D}^{(2)}\rightarrow\mathbb{R}$, where

\begin{align*} \mathcal{L}(q^{i},v^{A},\dot{v}^{A})=&\textcolor{blue}{\lambda_0}\mathcal{C}\left(q^{i},
v^{A},\dot{v}^{A}+\Gamma_{CB}^{A}v^{B}v^{C}+(\mathcal{G}^{\mathcal{D}})^{AB}\rho_{B}^{i}(q)\frac{\partial
V}{\partial q^{i}}\right)
\end{align*} \textcolor{blue}{where $\lambda_0\geq 0$.}

Then, the Lagrangian function $\mathcal{L}:\mathcal{D}^{(2)}\to\mathbb{R}$ is given by 

\begin{align*}
\mathcal{L}(q^{i},v^{A},\dot{v}^{C})&=\frac{\textcolor{blue}{\lambda_0}}{2}\left(||\gamma(t)- \gamma_{r}(t)||^2 + \epsilon ||u^A||^2\right)=\frac{1}{2}\left(\vphantom{\frac{\partial
V}{\partial q^{i}}\Big{|}\Big{|}^{2}}||q^{i}-q_r^{i}||^{2}+||v^{A}-v_r^{A}||^{2}\right.\\&\quad\left.+\epsilon||\Gamma_{CB}^{A}v^{C}v^{B}+\dot{v}^{A}+(\mathcal{G}^{\mathcal{D}})^{AB}\rho_{B}^{i}(q)\frac{\partial
V}{\partial q^{i}}\Big{|}\Big{|}^{2}\right)
\end{align*}

To derive the optimality conditions for the optimal tracking problem determined by $\mathcal{L}$ we use
standard variational calculus for systems with constraints by defining
the augmented Lagrangian $\widetilde{\mathcal{L}}
=\mathcal{L}-\lambda_{i}(\dot{q}^{i}-\rho_{A}^{i}(q)v^{A}).$
Therefore, the optimality conditions are given by the second-order Euler-Lagrange equations for $\widetilde{\mathcal{L}}$ (see \cite{Bl}, \cite{BlCoGuMdD}, \cite{leo}, \cite{leo2}) given by 
\begin{equation}\label{eqELmulti}
\dot{\lambda}_{i}=\frac{\partial\mathcal{L}}{\partial q^{i}}+\lambda_{j}\frac{\partial\rho_{A}^{j}}{\partial q^{i}}v^{A},\, \dot{q}^{i}=\rho_{A}^{i}(q)v^{A},\,
\frac{d}{dt}\left(\frac{\partial \mathcal{L}}{\partial\dot{v}^{A}}\right)=\frac{\partial \mathcal{L}}{\partial v^{A}}+\rho_{A}^{i}(q)\lambda_{i},\,
\end{equation} 

Observe that these equations arise from a constrained variational problem and the nonholonomic behavior is locally represented by the coordinates $(q^{i},v^{A})$ given by taking an adapted basis of vector fields in the nonholonomic distribution $\mathcal{D}$. The constraint enforced by the Lagrange multiplier \textcolor{blue}{$\lambda_i$} comes from the constraint arising from submanifold $\mathcal{D}^{(2)}$ \textcolor{blue}{and 
the solutions of the optimal control problem are the critical points of the functional $${\widetilde{J}}(\gamma, q,v,\dot{v}, \lambda,  \lambda_T)=\omega \Phi(T, \gamma(T))+\lambda_T r(\gamma(T), \gamma_r(T)) +\int_0^T\left[ 
\mathcal{L}-\lambda_{i}(\dot{q}^{i}-\rho_{A}^{i}(q)v^{A})\right]\,dt,$$
with $\omega>0$, $\gamma(0)\in {\mathcal D}$, $\lambda_T\in\mathbb{R}$ and $\gamma_r:[0, T]\rightarrow {\mathcal D}$ given. }

The optimal control problem for the nonholonomic system
given by $(\mathcal{D}^{(2)},\mathcal{L})$ with 
$\mathcal{L}:\mathcal{D}^{(2)}\ra\R$ is called \textit{regular} if and only if the
matrix
$\displaystyle{\left(\frac{\partial^{2}\mathcal{L}}{\partial\dot{v}^{A}\partial\dot{v}^{B}}\right)}$
is non singular (see \cite{BlCoGuMdD}, \cite{leo2}). For the proposed optimal trajectory tracking problem the system is always regular as long as $\epsilon\neq 0$. Note that our result coincides with the observation given in  \cite{OTT} Section $3.2$, and our Remark \ref{rksingular}, about when this class of optimal control problem becomes singular.

\begin{remark}The regularity condition is necessary to show the equivalence between the optimality conditions obtained by the variational approach and the ones obtained by employing the PMP as it was shown in \cite{BlCoGuMdD} (see Section  $4$ is \cite{BlCoGuMdD}) by using techniques of symplectic geometry. Therefore, since the optimal tracking problem  for the nonholonomic system
given by $(\mathcal{D}^{(2)},\mathcal{L})$ is regular, both formalisms are equivalent.\end{remark}

\subsection{Example: Optimal trajectory tracking for the nonholonomic particle}
Consider the situation of Example \ref{example}.

     The cost function $\mathcal{C}:\mathcal{D}\times\mathcal{U}\to\mathbb{R}$ for the optimal trajectory tracking problem is given by \begin{align*}
\mathcal{C}(q,v, u)= &\frac{\textcolor{blue}{\lambda_0}}{2} \left(|x- x_r |^2 + |y- y_r |^2 +|z- z_r |^2 \right.\\
&\left.+ |v^1- v^{1}_r|^2 + |v^2- v^{2}_{r}|^2 + \epsilon ((u^1)^2+(u^2)^2 )\right), 
\end{align*} and the terminal cost is determined by the function \begin{align*}\textcolor{blue}{r(\gamma(T),\gamma_r(T))}=
&|x(T)- x_r(T)|^2 + |y(T)- y_r(T) |^2 +|z(T)- z_r(T) |^2\\
+& |v^1(T)- v^{1}_r(T)|^2 + |v^2(T)- v^{2}_{r}(T)|^2\end{align*}with $T\in\mathbb{R}^{+}$ fixed.

Denoting by $(x,y,z,v^{1},v^{2},\dot{v}^{1},\dot{v}^{2})$ induced coordinates on $\mathcal{D}^{(2)}$ determined by the basis of vector fields $Y_1, Y_2$ which span $\mathcal{D}$ (see Example \ref{example}), the cost function $\mathcal{C}$ induces the Lagrangian $\mathcal{L}:\mathcal{D}^{(2)}\to\mathbb{R}$ given by
\begin{align*}
\mathcal{L}(q,v, \dot{v})= &\frac{\textcolor{blue}{\lambda_0}}{2} \left(\vphantom{\frac{2yv^{1}v^{2}\dot{v}^{2}}{1+y^{2}}}|x- x_r |^2 + |y- y_r |^2 +|z- z_r |^2 + |v^1- v^{1}_r|^2 + |v^2- v^{2}_{r}|^2 \right.\\
&\left.+ \epsilon (\dot{v}^1)^2+\epsilon\left((\dot{v}^2)^{2}+\frac{y^{2}}{(1+y^{2})^{2}}(v^{1}v^{2})^{2}+\frac{2yv^{1}v^{2}\dot{v}^{2}}{1+y^{2}}\right)\right),
\end{align*}  with $q=(x,y,z)$, $v=(v^{1}, v^2)$ and $\dot{v}=(\dot{v}^{1}, \dot{v}^{2})$. 

The extended Lagrangian is given by $$\widetilde{\mathcal{L}}(q,v, \dot{v})
=\mathcal{L}(q,v, \dot{v})-\lambda_{1}(\dot{x} +y v^2)-\lambda_2
     (\dot{y}- v^1)-\lambda_3(\dot{z}- v^2).$$

Necessary conditions for optimality are given by the solutions of the following system of nonlinear equations:  
\begin{align*}
\dot{\lambda}_1=&-\textcolor{blue}{\lambda_0}(x-x_r),\,\, \dot{\lambda}_3=-\textcolor{blue}{\lambda_0}(z-z_r)\\
\dot{\lambda}_2=&\epsilon\textcolor{blue}{\lambda_0} v^{1}v^{2}(y^{2}-1)\left(\frac{\dot{v}^{2}}{(1+y^{2})^{2}}+\frac{(v^1v^2) y}{(1+y^{2})^{3}}\right)+\lambda_1v^{2}-\textcolor{blue}{\lambda_0}(y-y_r),\\
\textcolor{blue}{\lambda_0}\epsilon\ddot{v}^{1}=&\textcolor{blue}{\lambda_0}(v^1-v^1_r)+\lambda_2+\frac{\textcolor{blue}{\lambda_0}\epsilon yv^{2}\dot{v}^{2}}{(1+y^{2})}+\frac{\textcolor{blue}{\lambda_0}\epsilon v^{1}(yv^{2})^{2}}{(1+y^{2})^2},\\
\textcolor{blue}{\lambda_0}\epsilon\ddot{v}^{2}=&\textcolor{blue}{\lambda_0}(v^2-v^2_r)-\lambda_1y+\lambda_3+\frac{2\textcolor{blue}{\lambda_0}\epsilon yv^{1}}{1+y^{2}}\left(\frac{yv^{1}v^{2}}{1+y^{2}}+\dot{v}^{2}\right),
\end{align*} together with  the admissibility conditions $\dot{x}= -y v^2$,
     $\dot{y}= v^1$ and 
     $\dot{z}= v^2.$

\section{Construction of variational integrators}\label{sec5}

Variational integrators (see \cite{mawest} for details) are derived
from a discrete variational principle. These integrators  retain
some of the main geometric properties of the continuous systems,
such as symplecticity, momentum conservation (as long as the
symmetry survives the discretization procedure), and good (bounded) behavior
of the energy associated to the system. 
of these type of variational integrators.

A \textit{discrete Lagrangian} is a differentiable function
$L_d\colon Q \times Q\to \R$, which may be considered as an
approximation of the action integral defined by a continuous regular 
Lagrangian $L\colon TQ\to \R.$ That is, given a time step $h>0$
small enough,
\[
L_d(q_0, q_1)\approx \int^h_0 L(q(t), \dot{q}(t))\; dt,
\]
where $q(t)$ is the unique solution of the Euler-Lagrange equations
for $L$ with  boundary conditions $q(0)=q_0$ and $q(h)=q_1$.

We construct the grid $\{t_{k}=kh\mid k=0,\ldots,N\},$ with $Nh=T$
and define the discrete path space
$\mathcal{P}_{d}(Q):=\{q_{d}:\{t_{k}\}_{k=0}^{N}\ra Q\}.$ We
identify a discrete trajectory $q_{d}\in\mathcal{P}_{d}(Q)$ with its
image $q_{d}=\{q_{k}\}_{k=0}^{N}$, where $q_{k}:=q_{d}(t_{k})$. The
discrete action $\mathcal{A}_{d}:\mathcal{P}_{d}(Q)\ra\R$ for this
sequence is calculated by summing the discrete Lagrangian on each
adjacent pair and is defined by
\begin{equation}\label{acciondiscreta}
\mathcal{A}_d(q_{d}) = \mathcal{A}_d(q_0,...,q_N) :=\sum_{k=0}^{N-1}L_d(q_k,q_{k+1}).
\end{equation}
We would like to point out that the discrete path space is
isomorphic to the smooth product manifold which consists of $N+1$
copies of $Q$. The discrete action inherits the smoothness of the
discrete Lagrangian and the tangent space
$T_{q_{d}}\mathcal{P}_{d}(Q)$ at $q_{d}$ is the set of maps
$v_{q_{d}}:\{t_{k}\}_{k=0}^{N}\ra TQ$ such that $\tau_{Q}\circ
v_{q_{d}}=q_{d}$ which will be denoted by
$v_{q_{d}}=\{(q_{k},v_{k})\}_{k=0}^{N},$ where $\tau_{Q} : TQ
\rightarrow Q$ is the canonical projection.

For any product manifold $Q_1\times Q_2,$
$T^{*}_{(q_1,q_2)}(Q_1\times Q_2)\simeq T^{*}_{q_1}Q_1\oplus
T^{*}_{q_2}Q_2,$ for $q_1\in Q_1$ and $q_2\in Q_2$ where $T^{*}Q$
denotes the cotangent bundle of a differentiable manifold $Q.$
Therefore, any covector $\alpha\in T^{*}_{(q_1,q_2)}(Q_1\times Q_2)$
admits an unique decomposition $\alpha=\alpha_1+\alpha_2$ where
$\alpha_i\in T^{*}_{q_i}Q_i,$ for $i=1,2.$ Thus, given a discrete
Lagrangian $L_d$ we have the following decomposition
  $$dL_{d}(q_0,q_1)=D_{1}L_d(q_0,q_1)+D_{2}L_d(q_0,q_1),$$
where $D_{1}L_d(q_0,q_1)\in T^*_{q_0}Q$ and $D_{2}L_d(q_0,q_1)\in
T^*_{q_1}Q$.

The discrete variational principle, states that
the solutions of the discrete system determined by $L_d$ must
extremize the action sum given fixed points $q_0$ and $q_N.$
Extremizing $\mathcal{A}_d$ over $q_k$ with $1\leq k\leq N-1,$ we
obtain the following system of difference equations
\begin{equation}\label{discreteeq}
 D_1L_d( q_k, q_{k+1})+D_2L_d( q_{k-1}, q_{k})=0.
\end{equation}

These equations are usually called \textit{the discrete Euler-Lagrange
equations}. Given a solution $\{q_{k}^{*}\}_{k\in\mathbb{N}}$ of
eq.\eqref{discreteeq} and assuming the regularity hypothesis (the
matrix $(D_{12}L_d(q_k, q_{k+1}))$ is regular), it is possible to
define implicitly a (local) discrete flow $
\Upsilon_{L_d}\colon\mathcal{U}_{k}\subset Q\times Q\to Q\times Q$
by $\Upsilon_{L_d}(q_{k-1}, q_k)=(q_k, q_{k+1})$ from
(\ref{discreteeq}), where $\mathcal{U}_{k}$ is a neighborhood of the
point $(q_{k-1}^{*},q_{k}^{*})$.

In order to construct structure-preserving variational integrators for nonholonomic mechanical control systems, one starts by considering the Lagrangian function $\mathcal{L}:\mathcal{D}^{(2)}\to\mathbb{R}$, where $\mathcal{D}^{(2)}$ is the submanifold of $T\mathcal{D}$. For simplicity in our computations,  from now on, we assume $Q$ is a real finite dimensional vector space. The tangent bundle of $\mathcal{D}$ can be discretized as $\mathcal{D}\times\mathcal{D}$. We define the submanifold $\mathcal{D}_{d}^{(2)}$ of $\mathcal{D}\times\mathcal{D}$ as 
\begin{equation*}
\mathcal{D}_{d}^{(2)}=\left\{(q_0^{i},v_{0}^{A},q_{1}^{i},v_{1}^{A})\in\mathcal{D}\times\mathcal{D}\mid\frac{q_1^{i}-q_0^{i}}{h}=\rho_{A}^{i}\left(\frac{q_{0}^{i}+q_{1}^{i}}{2}\right)\left(\frac{v_{0}^{A}+v_{1}^{A}}{2}\right)\right\},\end{equation*}
representing the discretization of $\mathcal{D}^{(2)}\subset T\mathcal{D}$. We assume that $Q$ is a vector space everywhere. 

One then discretizes the Lagrangian $\mathcal{L}:\mathcal{D}^{(2)}\to\mathbb{R}$ (we only discuss the mid-point rule here) as $\mathcal{L}_{d}:\mathcal{D}_{d}^{(2)}\to\mathbb{R}$,
\begin{equation}\label{ltildad}
\mathcal{L}_{d}(q_{k}^{i},v^{A}_{k},q^{i}_{k+1},v^{A}_{k+1}) = h\mathcal{L}(q_{k+1/2}^{i},v^{A}_{k+1/2},v^{A}_{k,k+1}),\end{equation}
where $(q_{k}^{i},v^{A}_{k},q^{i}_{k+1},v^{A}_{k+1})\in\mathcal{D}_{d}^{(2)}$ and where we are using the notation $z_{k+1/2}=\frac{1}{2}(z_k+z_{k+1})$ and $z_{k,k+1}=\frac{1}{h}(z_{k+1}-z_k)$.

Note that the discretization \eqref{ltildad} is carried out after writing the continuous-time Lagrangian $\mathcal{L}$ as a function of $(q^{i},v^{A},\dot{v}^{A})$. 

The variational integrator for the optimal control problem of the nonholonomic system is determined by minimizing the discrete action sum  $$\displaystyle{\mathcal{A}_{d}(\{q_{k}\}_{k=0}^{N-1})=\sum_{k=0}^{N-1}\mathcal{L}_{d}(q_{k}^{i},v^{A}_{k},q^{i}_{k+1},v^{A}_{k+1})}$$ over the path $(q_1,\ldots,q_{N-1},v_1,\ldots,v_{N-1})$ given fixed initial and final points $q_0,v_0$ and $q_N,v_N$, respectively, and subject to the discrete constraint functions $\Psi_{d}^{j}:\mathcal{D}_{d}^{(2)}\to\mathbb{R}$ with $j=1,\ldots,n=\dim(Q)$ given by $$\Psi_{d}^{j}(q_k^{i},v_{k}^{A},q_{k+1}^{i},v_{k+1}^{A})=q_{k,k+1}^{i}-\rho_{A}^{j}(q_{k+1/2}^{i})(v_{k,k+1}^{A}).$$ 
By considering the extended discrete action sum $$\displaystyle{\widetilde{\mathcal{A}}_d(\{q_{k}\}_{k=0}^{N})=\mathcal{A}_d(\{q_{k}\}_{k=0}^{N})+\sum_{k=1}^{N-1}(\lambda_{j}^{k})^{T}\Psi_{d}^{j}(q_{k}^{i},v^{A}_{k},q^{i}_{k+1},v^{A}_{k+1})},$$ where $\lambda_{j}^{k}=(\lambda_1^{k},\ldots,\lambda_n^{k})\in\mathbb{R}^{n}$ are the Lagrange multipliers. By extremizing the extended discrete action sum, with respect to variations $\delta q^{i}_{k}$, $\delta v_{k}^{A}$  and $\delta\lambda_{j}^{k}$, given fixed initial and final points $q_0,q_N,v_0,v_N$, satisfying the constraints, and using discrete integration by parts, leads to the following discrete Euler-Lagrange equations:
\begin{align*}
0=&D_{1}\mathcal{L}_{d}(q_k^{i},v_k^{A},q_{k+1}^{i},v_{k+1}^{A})+D_3\mathcal{L}_{d}(q_{k-1}^{i},v_{k-1}^{A},q_{k}^{i},v_{k}^{A})\\&+
\lambda^{k}_{j}D_{1}\Psi_{d}^{j}(q_k^{i},v_k^{A},q_{k+1}^{i},v_{k+1}^{A})+\lambda_{j}^{k-1}D_3\Psi_{d}^{j}(q_{k-1}^{i},v_{k-1}^{A},q_{k}^{i},v_{k}^{A}),\\
0=&D_{2}\mathcal{L}_{d}(q_k^{i},v_k^{A},q_{k+1}^{i},v_{k+1}^{A})+D_4\mathcal{L}_{d}(q_{k-1}^{i},v_{k-1}^{A},q_{k}^{i},v_{k}^{A})\\&+
\lambda^{k}_{j}D_{2}\Psi_{d}^{j}(q_k^{i},v_k^{A},q_{k+1}^{i},v_{k+1}^{A})+\lambda_{j}^{k-1}D_4\Psi_{d}^{j}(q_{k-1},v_{k-1},q_{k},v_{k}),\\
0=&\Psi_{d}^{j}(q_k^{i},v_k^{A},q_{k+1}^{i},v_{k+1}^{A}),
\end{align*}
for $k=1,\ldots,N-1$ and $j=1,\ldots,n$ and where $D_{i}$ represents the derivative  with respect to the $i^{th}$ argument. \textcolor{blue}{Note that initial conditions must belong to $\mathcal{D}$ and, $(q_N, v_N)=\gamma_r(Nh)$ (which is equivalent to impose that the constraint $r$ holds in discrete time) and fix $(q_N,v_N)$ to $\Phi(T,\gamma(T))$ if we consider the final cost.}

If the matrix \[\mathcal{M}=\left( \begin{array}{ccc}
D_{13}\mathcal{L}_d &D_{14}\mathcal{L}_d  & D_{13}\Psi_d^{j}\\
D_{23}\mathcal{L}_d &D_{24}\mathcal{L}_d  & D_{14}\Psi_d^{j}\\
D_{23}\Psi_d^{j} &D_{24}\Psi_d^{j}  & 0
 \end{array} \right)\]is non singular, the condition for local solvability of the constrained system is fulfilled and by the implicit function theorem the last set of equations determines an implicit local flow map, giving rise to the update map $\Upsilon:\mathcal{D}_{d}^{(2)}\times\R^{n}\to\mathcal{D}_{d}^{(2)}\times\R^{n}$
$$ \Upsilon(q_{k-1}^{i},v_{k-1}^{A},q_{k}^{i},v_{k}^{A},\lambda^{k-1})=(q_k^{i},v_k^{A},q_{k+1}^{i},v_{k+1}^{A},\lambda^{k}).$$

\subsection{Example: the Chaplygin sleigh} Consider the Chaplyigin sleigh of Example \ref{chaplygin} but subject to input controls. As we saw in Example \ref{chch2} the controlled Euler-Lagrange
equations are given by 
\begin{equation*}
\dot{v}^1=-\frac{a\sqrt{m}}{J+ma^2}v^1v^2+u_1,\quad\dot{v}^2=\frac{a\sqrt{m}}{J+ma^2}(v^1)^2+u_2. \end{equation*}together with the admissibility conditions
\begin{equation}\label{eqadm3}
\dot{x}_1=\frac{\cos\theta}{\sqrt{m}}v^2,\,\dot{x}_2=\frac{\sin\theta}{\sqrt{m}}v^{1},\,\dot{\theta}=\frac{1}{\sqrt{J+ma^2}}v^{1}
\end{equation}and the nonholonomic constraint $v^3=0$.

Here, $\mathcal{D}^{(2)}$ is defined by $(x_1,x_2,\theta,v^1,v^2,\dot{x}_1,\dot{x}_2,\dot{\theta},\dot{v}^1,\dot{v}^2)\in T\mathcal{D}$, satisfying \eqref{eqadm3}. Then the optimal control problem consists of finding an admissible curve
satisfying the previous equations given boundary conditions on
$\mathcal{D}$ and minimizing the functional
\begin{align*}
\mathcal{J}(x_1,x_2,\theta,v^1,v^2,u_1,u_2)=&\int_{0}^{T}\frac{\textcolor{blue}{\lambda_0}}{2}(|x_1- (x_1)_r |^2 + |x_2- (x_2)_r |^2 +|\theta- \theta_r |^2\\& +|v^1- v^{1}_r|^2 + |v^2- v^{2}_{r}|^2)+\frac{\textcolor{blue}{\lambda_0}\epsilon}{2}\left(u_{1}^{2}+u_{2}^{2}
\right)\,dt\, \end{align*} for the cost function $\mathcal{C}:\mathcal{D}\times \mathcal{U}\ra\R$
given by
\begin{align*}\mathcal{C}(x_1,x_2,\theta,v^1,v^2,u_1,u_2)=&\frac{\textcolor{blue}{\lambda_0}\epsilon}{2}(u_1^2+u_2^2)+\frac{\textcolor{blue}{\lambda_0}}{2}(|x_1- (x_1)_r |^2 + |x_2- (x_2)_r |^2 \\&+|\theta- \theta_r |^2 + |v^1- v^{1}_r|^2 + |v^2- v^{2}_{r}|^2 ),\end{align*}where $\gamma(t)=(x_1(t),x_2(t),\theta(t),v^1(t),v^2(t))$ and also we must to take care that \textcolor{blue}{$\theta\in[0,2\pi)$}.

The optimal control problem is equivalent to solving the
constrained variational problem determined by
$\mathcal{L}:\mathcal{D}^{(2)}\ra\R,$ where

\begin{align}
\mathcal{C}(x_1,x_2,\theta,v^1,v^2,\textcolor{blue}{\dot{v}^{1}}, \textcolor{blue}{\dot{v}^{2}})= &\frac{\textcolor{blue}{\lambda_0}}{2} \left(|x_1- (x_1)_r |^2 + |x_2- (x_2)_r |^2 +|\theta- \theta_r |^2 \right.\\\label{costchaplygin2}
&\left.+ |v^1- v^{1}_r|^2 + |v^2- v^{2}_{r}|^2 \right)+ \textcolor{blue}{\lambda_0}\epsilon(\dot{v}^{1}+\eta v^1v^2)^2\\&+\textcolor{blue}{\lambda_0}\epsilon(\dot{v}^2-\eta(v^1)^2)^2 .\nonumber
\end{align} where $\displaystyle{\eta=\frac{a\sqrt{m}}{J+ma^2}}$. \textcolor{blue}{We also introduce the discrete version of constraint constraint $r(\gamma(T),\gamma_r(T))=0$ where
\begin{align*}
r_d(\gamma_{d,N},(\gamma_{r})_{d,N})=
&|x_{1,N}- (x_1)_{r,N} |^2 + |x_{2,N} - (x_2)_{r,N} |^2 +| \theta_{N} - \theta_{r,N} |^2\\
&+ |v^1_{N}- (v^{1})_{r,N}|^2 + |v^2_{N}- (v^{2})_{r,N}|^2\,,
\end{align*}
where $(\gamma_{r})_{d}$ denotes a discrete reference trajectory. This can be, for instance, an uncontrolled instance of the same system.}

 Consider the  extended Lagrangian
\begin{small}$$\widetilde{\mathcal{L}}=\displaystyle{\mathcal{L}+\lambda_{1}\left(\dot{x}_1-\frac{\cos\theta}{\sqrt{m}}v^2\right)+\lambda_{2}\left(\dot{x}_2-\frac{\sin\theta}{\sqrt{m}}v^2\right)+\lambda_{3}\left(\dot{\theta}-\frac{1}{\sqrt{J+ma^2}}v^{1}\right)}$$\end{small}
with \begin{align*}\mathcal{L}=&\frac{\textcolor{blue}{\lambda_0}}{2} \left(|x_1- (x_1)_r |^2 + |x_2- (x_2)_r |^2 +|\theta- \theta_r |^2+ |v^1- v^{1}_r|^2 + |v^2- v^{2}_{r}|^2 \right.\\&\left.+ \epsilon(\dot{v}^{1}+\eta v^1v^2)^2+\epsilon(\dot{v}^2-\eta(v^1)^2)^2 \right)\end{align*} where $\displaystyle{\eta=\frac{a\sqrt{m}}{J+ma^2}}$.

The optimality conditions are then given by 
\begin{align*}
\dot{\lambda}_{1}=&\textcolor{blue}{\lambda_0}(x_1-(x_1)_r),\, \dot{\lambda}_{2}=\textcolor{blue}{\lambda_0}(x_2-(x_2)_r),\,
\dot{\lambda}_{3}=\textcolor{blue}{\lambda_0}(\theta-\theta_r)+\lambda_1\frac{\sin\theta}{\sqrt{m}}v^2-\lambda_2\frac{\cos\theta}{\sqrt{m}}v^2,\\
\ddot{v}^1=&(\dot{v}^{1}+\eta v^1v^2)v^2\eta-2\eta v^1(\dot{v}^{2}-\eta(v^1)^2)+\frac{(v^1-v^1_r)}{\textcolor{blue}{\lambda_0}\epsilon}-\frac{\lambda_2\sin\theta}{\textcolor{blue}{\lambda_0}\epsilon\sqrt{m}}-\frac{\lambda_3}{\textcolor{blue}{\lambda_0}\epsilon\sqrt{J+ma^2}}\\&-\eta\dot{v}^{1}v^{2}-\eta v^1\dot{v}^2,\\
\ddot{v}_{2}=&2\eta v^1\dot{v}^1+(\dot{v}^1+\eta v^1v^2)\eta v^1-\frac{\lambda_1\cos\theta}{\textcolor{blue}{\lambda_0}\epsilon\sqrt{m}}+\frac{(v^1-v^2_r)}{\textcolor{blue}{\lambda_0}\epsilon},
\end{align*} together with the admissibility conditions \eqref{eqadm3}.

The variational integrator for the optimal control problem of the Chaplygin sleigh is constructed by the discretization of the Lagrangian \eqref{costchaplygin2} and the construction of the space $\mathcal{D}_{d}^{(2)}$ which determines the  discrete constraint.

Let $h\in\mathbb{R}^{+}$ be the time step.  To simulate solutions of the tracking problem we apply the mid-point rule to the cost function and  constraints,  for $h=0.1$ and $N=50$ intervals (and therefore $51$ nodes).

\begin{figure}[h!]\label{fig4}
	\centering
	\subfigure{\includegraphics[trim = {10mm 85mm 18mm 95mm}, scale=0.35, clip=true]{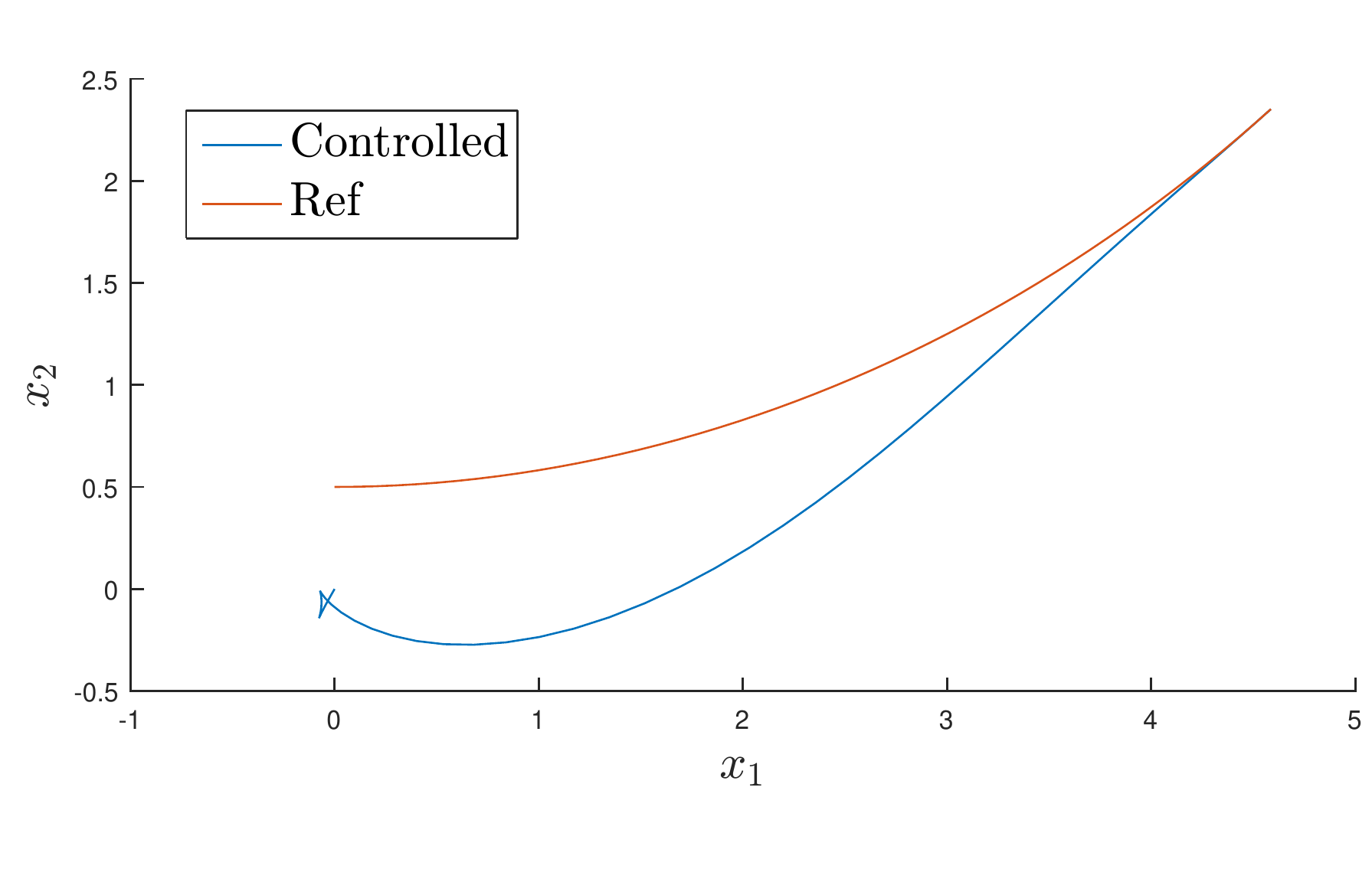}}
	\subfigure{\includegraphics[trim = {10mm 85mm 18mm 95mm}, scale=0.35, clip=true]{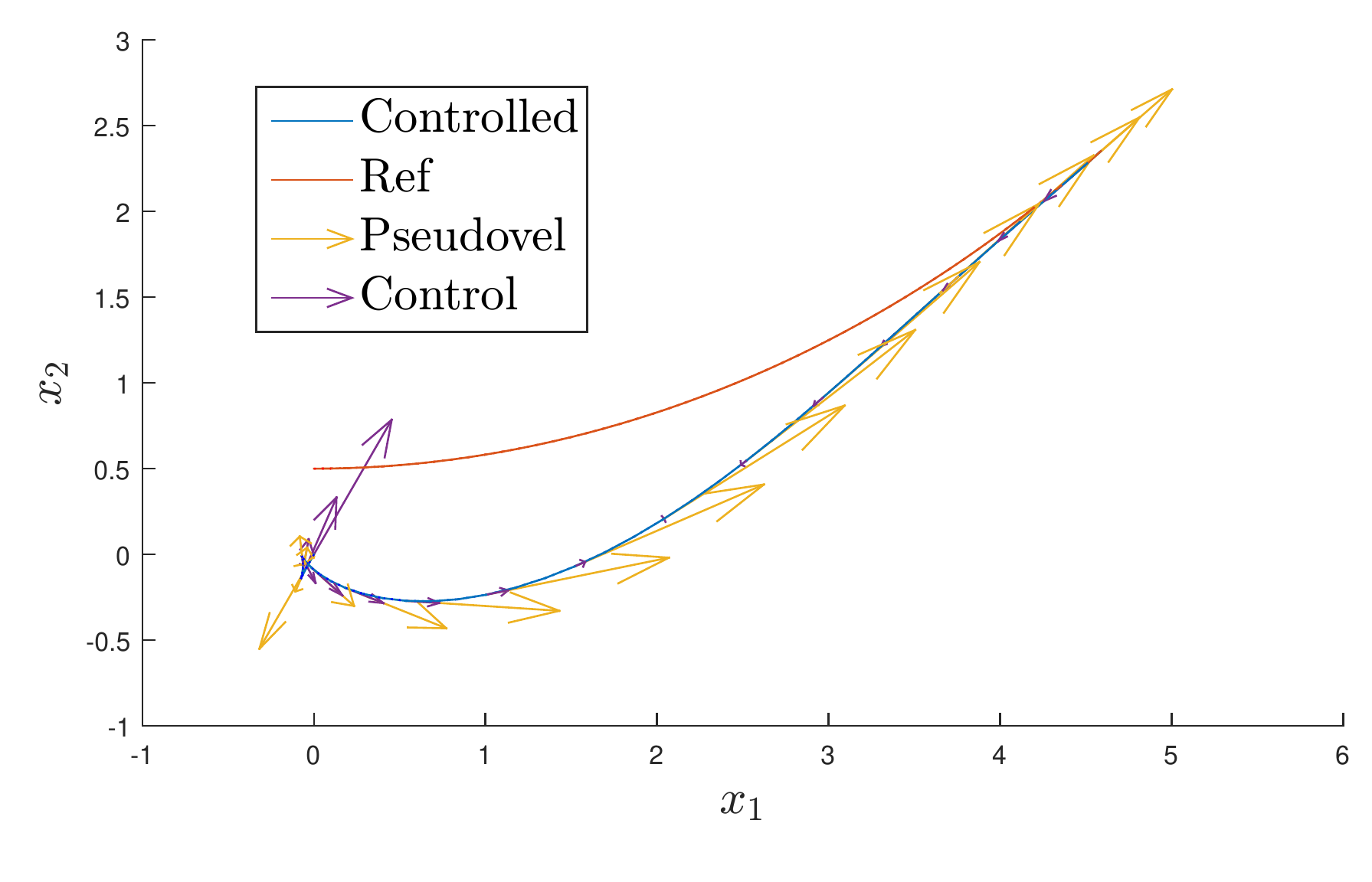}}
	\caption{Trajectories minimizing the cost function $\mathcal{J}$, evolving on $\mathcal{D}$ and tracking the reference trajectory $\gamma_r$ in time $T$. Note that the initial conditions of the controlled trajectory oblige it to stop its forward motion, back up and turn to correct its direction. Left: controlled trajectory in blue, reference trajectory $\gamma_r$ in red. Right: Superimposed quasivelocities in yellow and control vector field in purple.}
\end{figure}

\begin{figure}[h!]\label{fig5}
	\centering
	\subfigure{\includegraphics[trim = {10mm 75mm 18mm 85mm}, scale=0.35, clip=true]{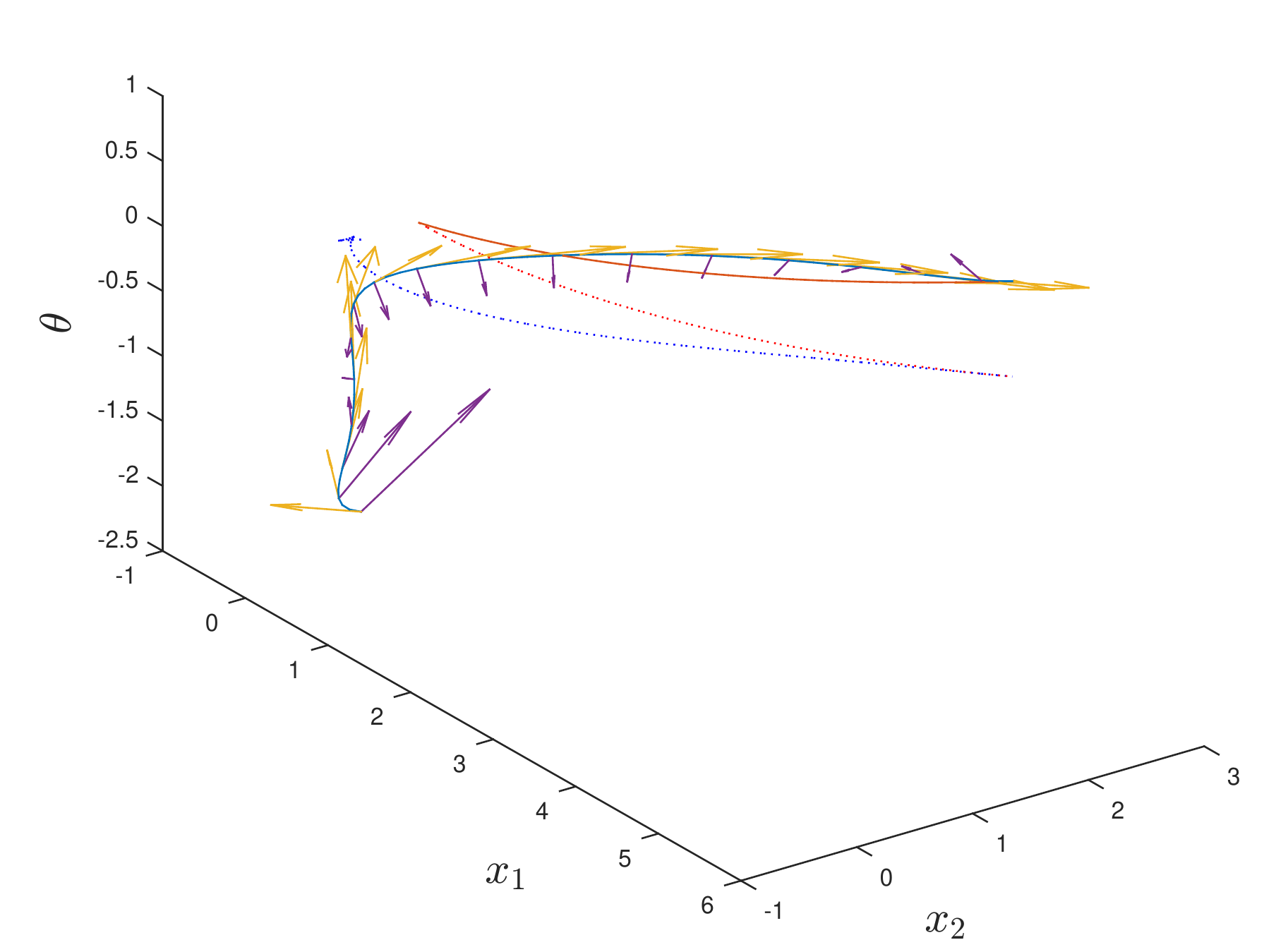}}
	\subfigure{\includegraphics[trim = {10mm 75mm 18mm 85mm}, scale=0.35, clip=true]{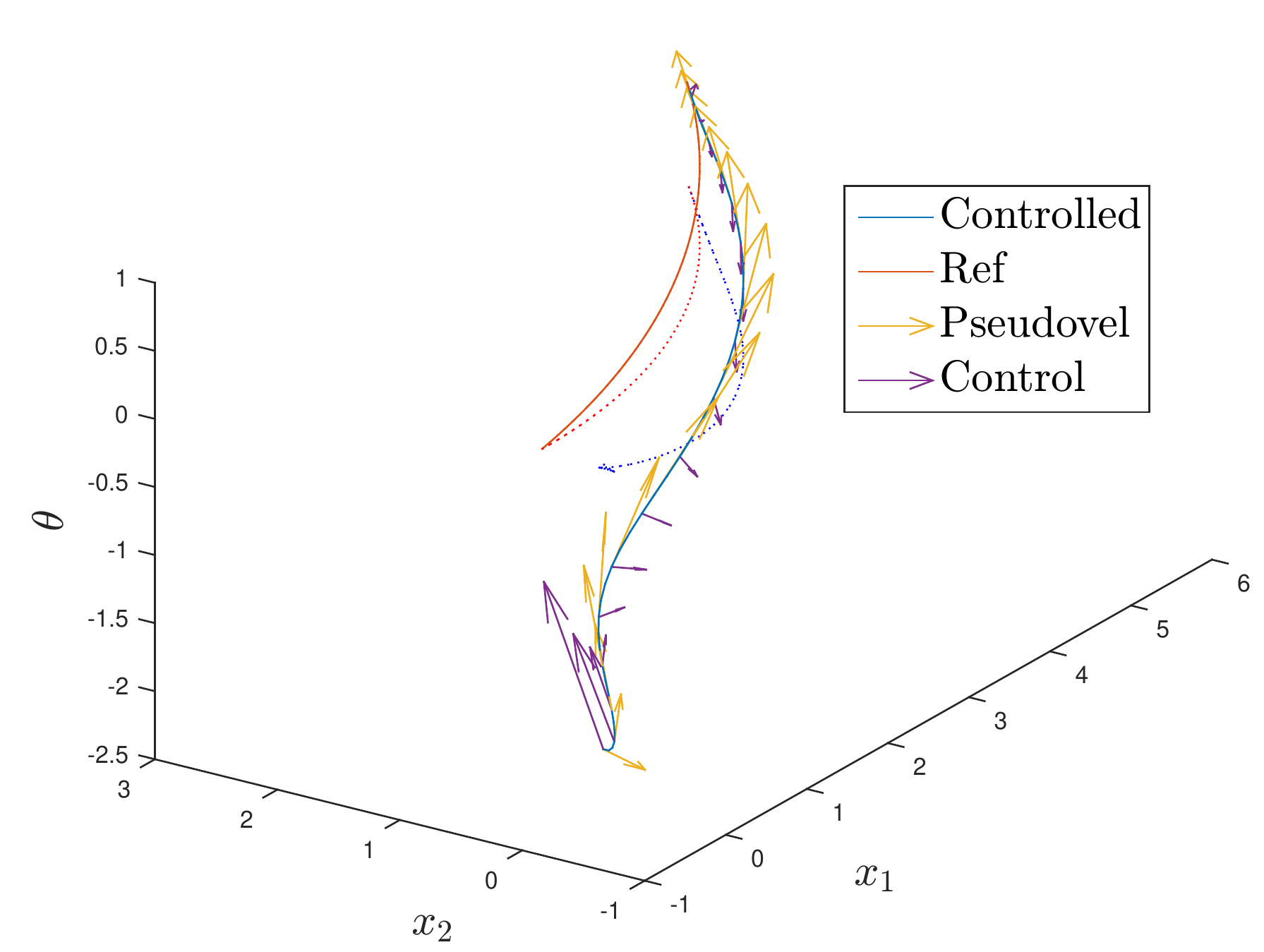}}
	\caption{Trajectories minimizing the cost function $\mathcal{J}$, evolving on $\mathcal{D}$ and tracking the reference trajectory $\gamma_r$ in time $T$. 3D representation with angle in vertical axis. The blue curve represents the controlled trajectory and the red curve the reference trajectory $\gamma_r$, the yellow vectors show the quasivelocities along the evolution of the curve and the purple vectors represent the control vector field. The dotted lines are the planar projection of the trajectories onto the $\theta = 0$ plane.}
\end{figure}

\begin{figure}[h!]\label{fig6}
	\centering
	\subfigure{\includegraphics[trim = {38mm 90mm 40mm 100mm}, scale=0.45, clip=true]{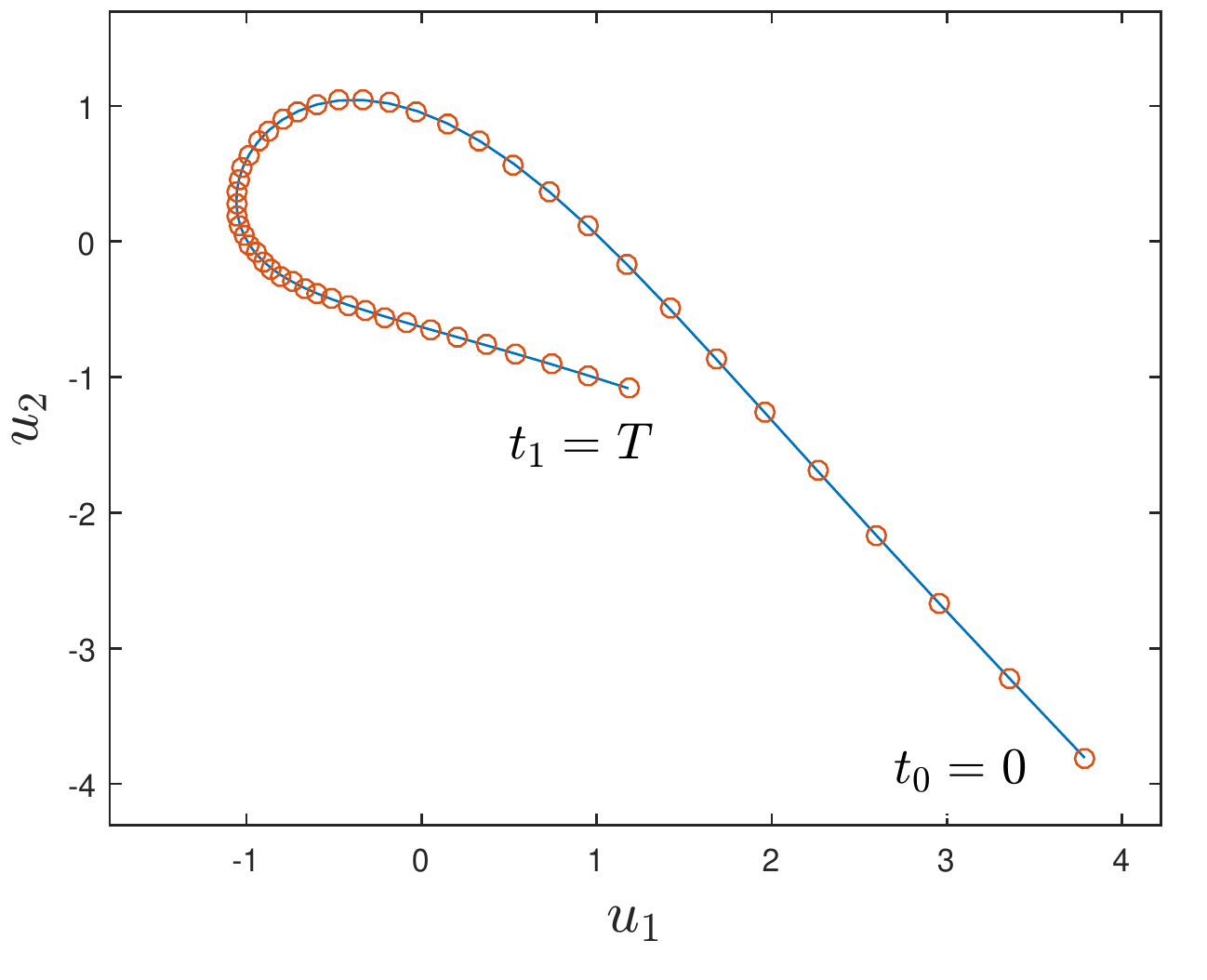}}\quad
	\subfigure{\includegraphics[trim = {38mm 90mm 40mm 100mm}, scale=0.45, clip=true]{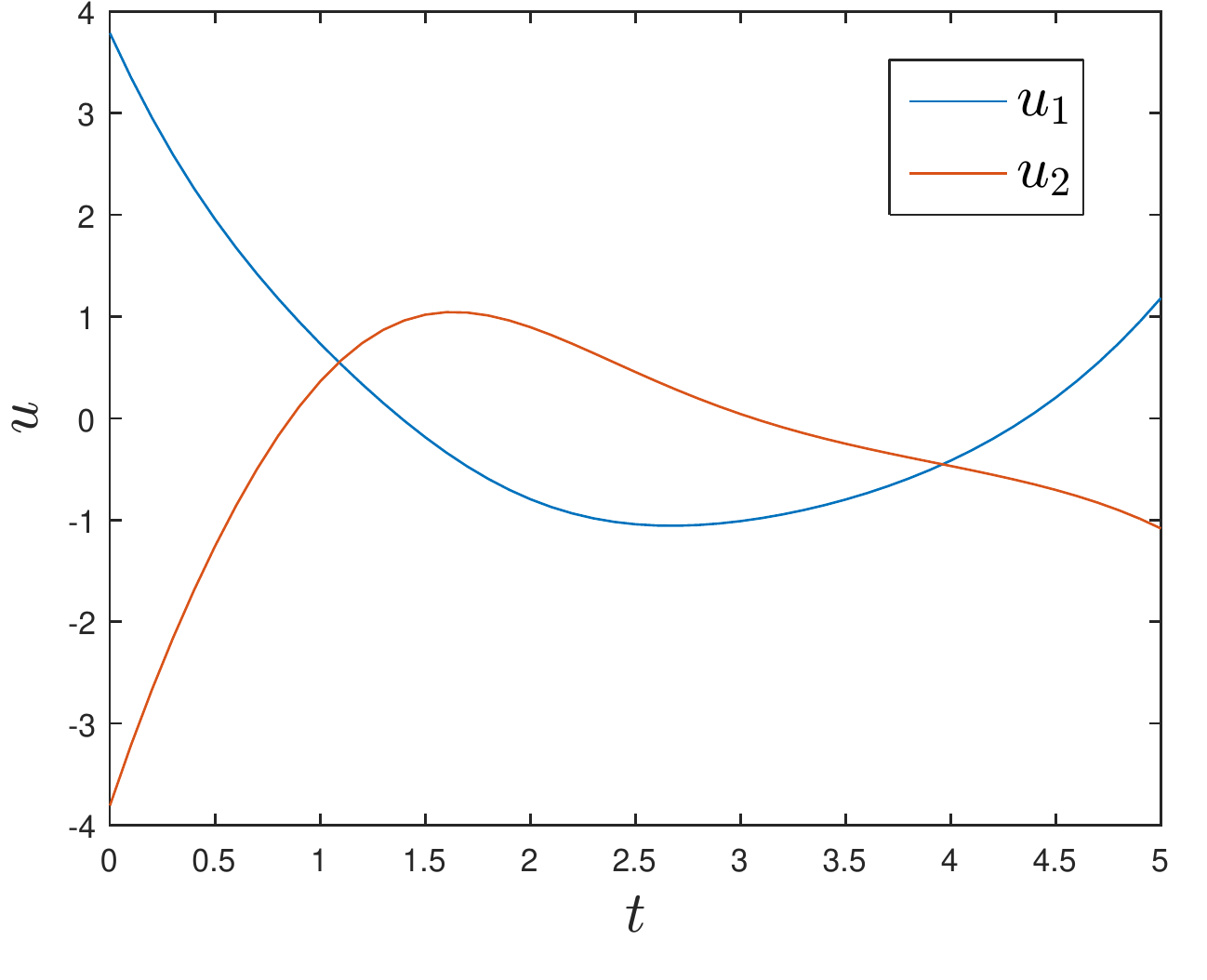}} 
	\caption{Control inputs minimizing the cost function $\mathcal{J}$, evolving on $\mathcal{D}$ and tracking the reference trajectory $\gamma_r$ in time $T$. Left: Representation of the control curve $(u_1, u_2)$, with red circles marking each time step. Right: Time evolution of the controls.}
\end{figure}

\begin{figure}[h]\label{fig7}
	\centering
	\subfigure{\includegraphics[trim = {38mm 90mm 40mm 100mm}, scale=0.45, clip=true]{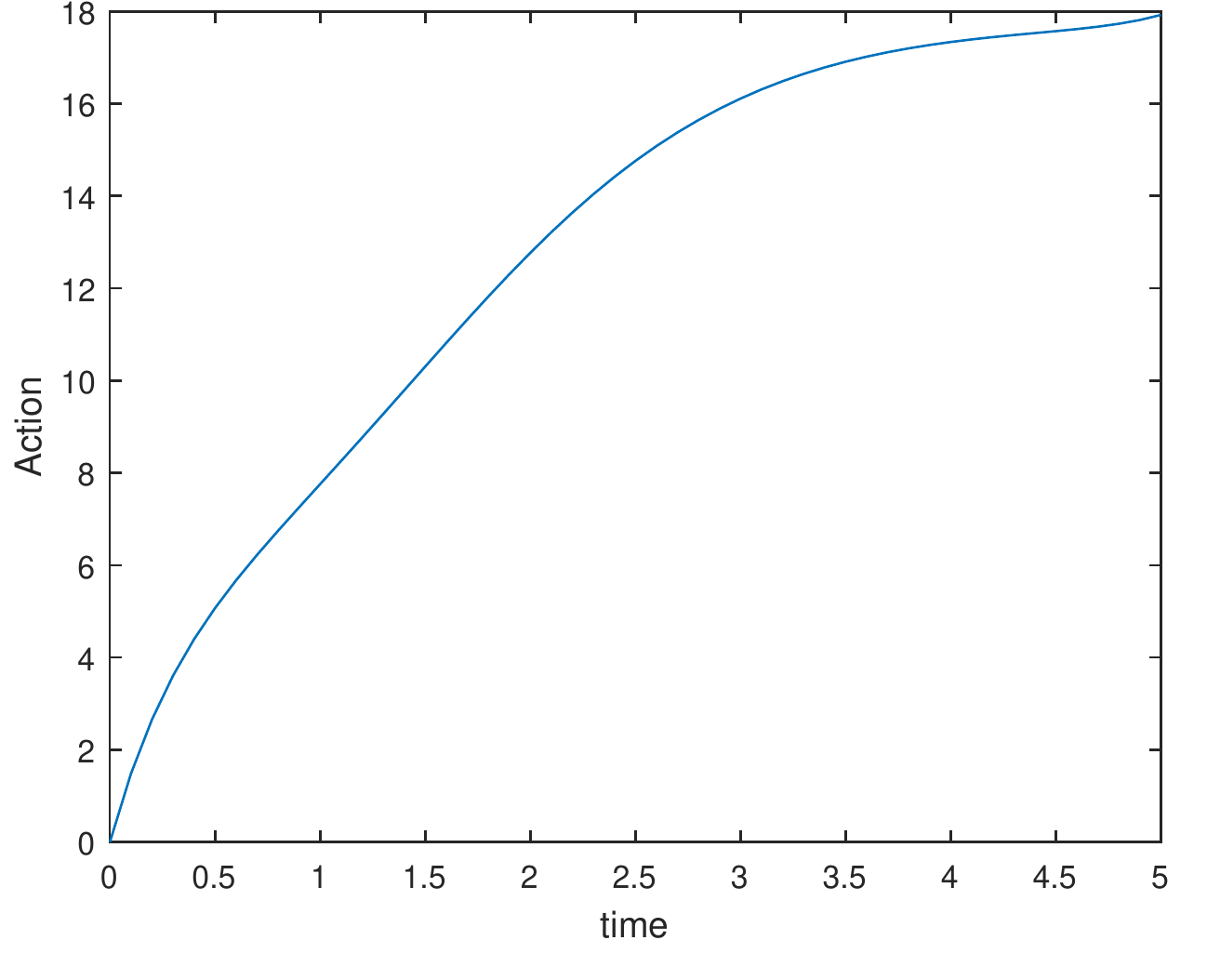}}\quad
	\subfigure{\includegraphics[trim = {38mm 90mm 40mm 100mm}, scale=0.45, clip=true]{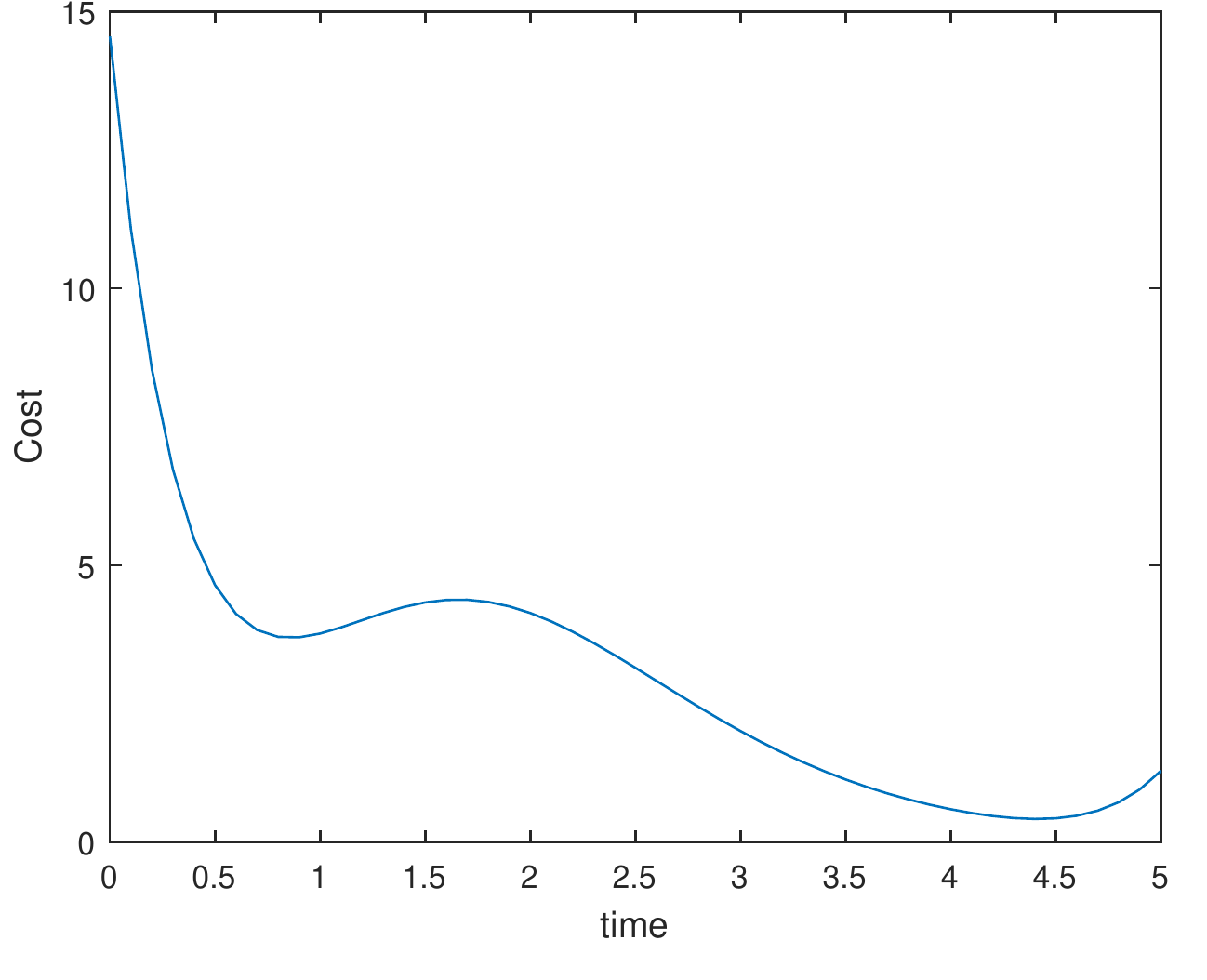}}
	\caption{Left: Time evolution of the action integral $\mathcal{J}$ using our variational integrator. Right: Time evolution of the cost function $\mathcal{C}$ using our variational integrator.}
\end{figure}

\begin{figure}[h]\label{fig8}
	\centering
	\includegraphics[trim = {15mm 85mm 20mm 90mm}, scale=0.45, clip=true]{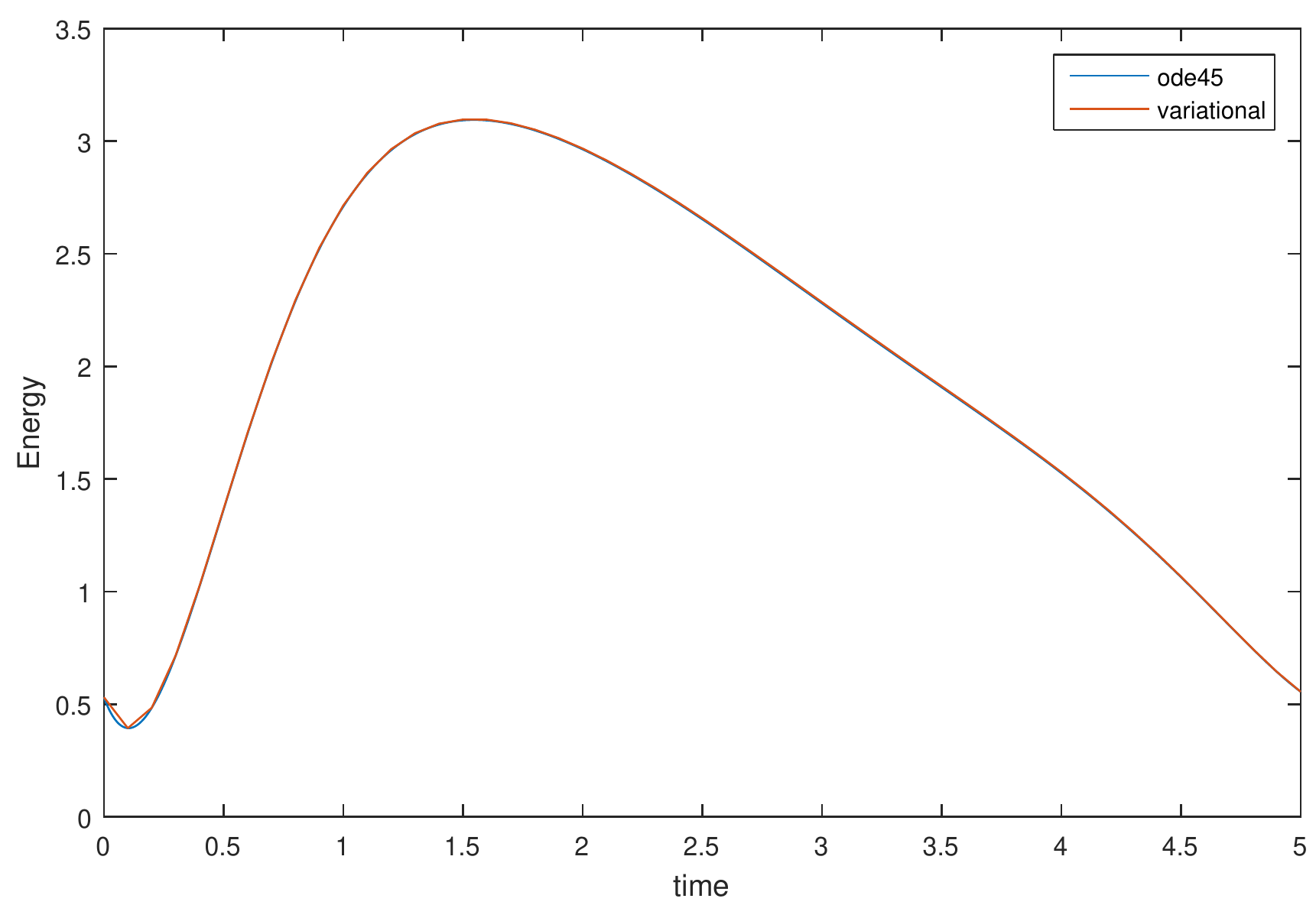}
	\caption{Time evolution comparison of the energy of the controlled sleigh. The blue line represents the one obtained via MATLAB's {\normalfont\texttt{ode45}} and the red line via our variational method. Note that our discretization is coarser (only 51 equidistant points) but it still manages to capture the behaviour remarkably well.}
\end{figure}

For the initial condition $\gamma(0)=(x^0_1,x^0_2,\theta_0,v_0^1,v_0^2)=(0,0,4\pi/3,1/4,1)$, \textcolor
{blue}{$\lambda_0=1$, $\lambda(T)=\lambda_T$ arbitrary for the shooting} \textcolor{blue}{and the reference trajectory $\gamma_r(t)$ is the uncontrolled trajectory of a Chaplygin sleigh with} $\gamma_r(\textcolor{blue}{0})=(x^{ref}_1,x^{ref}_2,\theta_{ref},v_{ref}^1,v_{ref}^2)=(0,1/2,0,1/3,1)$, $T=5$, $\epsilon=1$,  $m=1$, $J=4$ and $a=0.2$, $\lambda^0=\textcolor{blue}{\lambda(t_0)= 0}$, we exhibit the results in Figures \ref{fig4}, \ref{fig5}, \ref{fig6}, \ref{fig7} and \ref{fig8}.

The controlled generated by our trajectory planning to track the desires configurations have not been assessed in terms of their stability; we would,
therefore, like to find a method for incorporating the stability
of the nonholonomic system into our methodology. Similarly, it would be
of interest to study the cost of tracking them as a reference
trajectory. Finally, the  method proposed in this work can only guarantee
local optimality, and in our simulations the controlled Chaplygin sleigh displayed a multitude of local minima. Incorporating discrete mechanics into methods seeking the global optimum
of a cost functional, or bounds on it, remains an open task.

\section{Final Discussion} \label{sec6}

A class of nonlinear optimal control problems has
been identified to study tracking of trajectories for nonholonomic systems after detecting fundamental issues in the study of the  error dynamics applied to these problems. The nonlinear features arise directly from physical assumptions about constraints and Lagrangian dynamics on the motion of a mechanical system. The geometric framework introduced permits to study mechanical systems reduced by Lie group symmetries and multi-agent systems \cite{codi}, which will be further developed in an extension of this work, as well as variational interpolation problems \cite{interpolation}. We have
studied how to employ a shooting method and identify control issues for this class of systems and, we have derived new insights in this fundamental problem based on optimal control theory and tracking of trajectories.  The general approach described on this paper makes
substantial use of the geometric approach to nonlinear
control. However, the specific nonlinear control strategy
suggested is substantially different, both conceptually and
in detail, from the smooth nonlinear control strategies
most commonly studied in the literature.

 Minimizing the cost function while evolving on the constraint submanifold and remaining differentiable by solving a boundary value problem using a single shooting method is a difficult task and not very numerically stable, and this without considering time as an independent variable, which will only complicate things further.  In this work we consider tracking a trajectory as being synonymous with converging into it in a finite and prescribed time. Nevertheless, we believe that the optimality of the method may be improved by considering time as an additional degree of freedom, and setting the final time as a free and optimizable. This extension will be considered in a further publication. Next by analyzing the convergence to the reference trajectory by modifying the problem statement for a time horizon problem will be explored. The idea is to include an external dissipative force and study the problem by employing the dynamic programming principle and approximate the infinite time horizon problem with a the finite horizon  problem with terminal cost as in \cite{saccon}.

\textbf{Acknowledgments:} The authors wish to thank Prof. Ravi Banavar for fruitful comments about the preliminary version of this work. The authors
are indebted with the reviewers and editor for their recommendations that
helped to improve the quality, clarity and exposition of this
work. 







\begin{thebibliography}{99}


\bibitem{Bl} A. M. Bloch. \textit{Nonholonomic Mechanics and Control.} Interdisciplinary
Applied Mathematics Series, 24, Springer-Verlag, New York (2003).


\bibitem{bl1}A. M. Bloch. Stabilizability of nonholonomic control systems.
 Automatica, vol. 28, no. 2, pp. 431-435, 1992.
 
 \bibitem{BlCoGuMdD} A. Bloch, L. Colombo, R. Gupta and D. Mart\'in de Diego.  A Geometric Approach to the Optimal Control of Nonholonomic Mechanical Systems. \textit{Analysis and Geometry in Control Theory and its Applications}. INdAM series. Vol 12. 2015. 

\bibitem{interpolation} A. Bloch, M. Camarinha and L. J. Colombo. Dynamic interpolation for obstacle avoidance on Riemannian manifolds. International Journal of Control,
pages 1-22, doi:10.1080/00207179.2019.1603400. Preprint available at. arXiv:1809.03168
[math.OC].

\bibitem{bl2}A. M. Bloch and N. H. McClamroch. Control of mechanical
systems with classical nonholonomic constraints. Proc. IEEE
Conference Decision and Control, 1989, Tampa, FL, pp. 201-205.

\bibitem{bl3}A. Bloch, N. McClamroch, and M. Reyhanoglu. Controllability
and stabilizability properties of a nonholonomic control
system. Proc. IEEE Conference on Decision and Control, 1990,  1312-1314.



\bibitem{blochnonh}A. M Bloch, M. Reyhanoglu, and N. H. McClamroch. Control
and stabilization of nonholonomic dynamic systems. IEEE Transactions on
Automatic control, 37(11):1746-1757, 1992

\bibitem{BloZen} A. M. Bloch, J.E. Marsden and D. Zenkov. Quasivelocities and symmetries in non-holonomic
systems. Dynamical Systems, 24 (2), (2009),  187--222.

\bibitem{brocket} R. W. Brockett. Control theory and singular Riemannian geometry.
New Directions in Applied Mathematics, P. J. Hilton and
G. S. Young, Eds. New York Springer-Verlag, 1982. 
\bibitem{brocket2}R. W. Brocket. Asymptotic stability and feedback stabilization. Differential Geometric Control Theory. R. W. Brockett, R. S. Millman, and H. J.
Sussmann, Eds. Boston, MA: Birkhauser, 1983.
\bibitem{bullolewis} F. Bullo and A. Lewis. Geometric Control of Mechanical Systems: Modeling, Analysis, and Design for Simple Mechanical Control Systems. Texts in Applied Mathematics, Springer Verlag 2005.
\bibitem{celledoni} E. Celledoni, M. Farre Puiggali, E. Hoiseth, D. Martin de Diego. Energy-preserving
integrators applied to nonholonomic systems, Journal of Nonlinear Science, 29(4), 1523-1562, 2019.
\bibitem{leothesis} L. Colombo. Geometric and numerical methods for optimal control of mechanical systems. PhD thesis, Instituto de Ciencias Matem\'aticas, ICMAT (CSICUAM-UCM-UC3M), 2014.
\bibitem{leo}L. Colombo. A variational-geometric approach for the optimal control of nonholonomic systems. \textit{International Journal of Dynamics and Control}. Vol 6 (2), 652-662, 2018.
\bibitem{codi}  L.J. Colombo and D.V. Dimarogonas. Motion Feasibility Conditions for Multi-Agent Control Systems on Lie Groups. in IEEE Transactions on Control of Network Systems.
doi: 10.1109/TCNS.2019.2925264. Preprint available at arXiv preprint. arXiv:1808.04612, 2018.
\bibitem{leo2} L Colombo, R Gupta, A Bloch, DM de Diego. Variational discretization for optimal control problems of nonholonomic mechanical systems. Decision and Control (CDC), 2015 IEEE 54th Annual Conference on, 4047-4052.
\bibitem{Cobook}J. Cort\'es. Geometric control, and numerical aspects of nonholonomic systems. Lecture notes in Mathematics, Springer Verlag, 2002.
\bibitem{CoMa}
J. Cort\'es and E. Mart\'inez E. \textit{Mechanical control systems on Lie algebroids}. IMA J. Math. Control. Inf. 21, 457-492, 2004. 
\bibitem{Ha}H Hajieghrary, D Kularatne, M.A. Hsieh. Differential Geometric Approach to Trajectory Planning: Cooperative Transport by a Team of Autonomous Marine Vehicles. arXiv:1805.00959.


\bibitem{nijmeijer97} Z.-P. Jinag and H. Nijmeijer. Tracking control of mobile robots:
A case study in backstepping. Automatica, 33(7):1393-1399, 1997.

\bibitem{kodi} D. Koditschek. The application of total energy as a Lyapunov function
for mechanical control systems. Contemporary Math. 97-131,
1989
\bibitem{lewis} F. Lewis. \textit{Optimal control}. John Wiley $\&$ Sons, Inc, 1986.
\bibitem{OTT} J. L\"ober. \textit{Optimal trajectory tracking}. PhD. Thesis. TU Berlin, 2015.
\bibitem{mala}
\textcolor{blue}{
K. Malanowski. \textit{On normality of {L}agrange multipliers for state constrained
	optimal control problems}, Optimization. A Journal of Mathematical Programming and
	Operations Research,
52(1): 75-91, 2003.
}


\bibitem{mawest}
  J. E. Marsden and M. West.
  \textit{Discrete Mechanics and variational integrators.}
  Acta Numerica Vol.10 (2001), 357--514.
  

\bibitem{aradhana1}
A. Nayak and R. N. Banavar. 
On Almost-Global Tracking for a Certain Class of Simple Mechanical Systems.  IEEE Transactions on Automatic Control, Vol. 64 (1), 412-419, 2019.

\bibitem{aradhana2} 
A. Nayak, R. N. Banavar and D. H. S. Maithripala. Almost-global tracking for a rigid body with internal rotors. European Journal of Control
Vol 42, 59-66, 2018.
\bibitem{ACC} A. Nayak, R. T. Sato Mart\'in de Almagro, L. Colombo, D. Mart\'in de Diego. Optimal Trajectory Tracking of Nonholonomic Mechanical Systems: a geometric approach. 2019 American Control Conference (ACC), Philadelphia, PA, USA, 2019, pp. 1924-1929. doi: 10.23919/ACC.2019.8814647
Preprint available online at arXiv:1901.10374 [cs.SY]


\bibitem{saccon} A. Saccon, J. Hauser, A. P. Aguilar. Exploration of Kinematic Optimal Control on the Lie Group $SO (3)$. 8th IFAC Symposium on Nonlinear Control Systems. 1302-1307, 2010.


\bibitem{sanyal} A. Sanyal, N. Nordkvist, M. Chyba. An almost global tracking control scheme for maneuverable autonomous vehicles and its discretization.
IEEE Transactions on Automatic control, 56.2:457-462, 2011.

\bibitem{sina} S. Ober-Blobaum. Galerkin variational integrators and modified symplectic Runge--Kutta methods. IMA Journal of Numerical Analysis 37 (1), 375-406, 2017.
\bibitem{zh} A. Zuyev. Exponential Stabilization of Nonholonomic Systems by Means of Oscillating Controls. SIAM Journal on Control and Optimization, 2016, Vol. 54, No. 3: pp. 1678-1696.


\end{thebibliography}
\end{document}